\documentclass[12pt]{amsart}
\usepackage{amssymb}
\usepackage{epsf}
\usepackage{a4}

\numberwithin{equation}{section}

\newtheorem{thm}{Theorem}
\newtheorem{prop}[thm]{Proposition}
\newtheorem{lemma}[thm]{Lemma}
\newtheorem{cor}[thm]{Corollary}

\newtheorem{example}[thm]{Example}
\newtheorem{remark}[thm]{Remark}

\newenvironment{ex}{\begin{example}\rm}{\end{example}}
\newenvironment{rem}{\begin{remark}\rm}{\end{remark}}
\newcounter{FNC}[page]
\def\newfootnote#1{{\addtocounter{FNC}{1}$^\fnsymbol{FNC}$%
     \let\thefootnote\relax\footnotetext{$^\fnsymbol{FNC}$#1}}}
\newcommand{\G}{\mathbb{G}}
\newcommand{\R}{\mathbb{R}}

\renewcommand{\P}{\mathbb{P}}

\title[Lines Meeting a fixed line and Tangent to Two Spheres]{The Envelope of
       Lines Meeting a fixed line and Tangent to Two Spheres} 

\author{G\'abor Megyesi}
\address{Department of Mathematics\\
        UMIST\\
        P.O. Box 88\\
        Manchester, M60 1QD \\
        England}
\email{gmegyesi@ma.umist.ac.uk}
\urladdr{http://www.ma.umist.ac.uk/gm/}

\author{Frank Sottile}
\address{Department of Mathematics\\
         Texas A\&M University\\
         College Station\\
         TX \ 77843\\
         USA}
\email{sottile@math.tamu.edu}
\urladdr{http://www.math.tamu.edu/\~{}sottile}

\thanks{
Research of first author supported in part by EPSRC grant GR/S11381/01.
Research of second author supported in part by NSF grant DMS-0134860.}
\begin{document}

\begin{abstract}
 We study the set of lines that meet a fixed line and are tangent to two
 spheres and classify the configurations consisting of a single line
 and three spheres for which there are infinitely many lines tangent to the
 three spheres that also meet the given line.
 All such configurations are degenerate.
 The path to this result involves the interplay of some beautiful and
 intricate geometry of real surfaces in 3-space, complex projective algebraic 
 geometry, explicit computation and graphics. 
\end{abstract}

\maketitle

\section*{Introduction}

We determine the configurations of one line and three spheres for which there
are infinitely many common tangent lines to the spheres that also meet the
fixed line. 
Configurations of four lines having infinitely many common transversal lines
were described classically
and Theobald~\cite{Th02a} treated configurations of three lines and one sphere with 
infinitely many lines tangent to the sphere that also meet the fixed lines.
The case of two lines and two spheres was solved in~\cite{2l2s}. 

Besides the beautiful geometry encountered in this
study,\addtocounter{FNC}{1}\newfootnote{See: 
{\tt www.math.tamu.edu/\~{}sottile/pages/2l2s/} \ and \ 
 {\tt \~{}sottile/stories/3Spheres/}}
these questions were motivated by algorithmic problems
in computational geometry. 
As explained in~\cite{Th02a}, problems of this type occur 
when one is looking for a line or ray interacting 
(in the sense of ``intersecting'' or 
in the sense of ``not intersecting'') with a given set of three-dimensional
bodies, if the class of admissible bodies consists of polytopes
and spheres. 
Concrete instances include visibility computations with moving  
viewpoints~\cite{Th02b}, controlling a laser beam in
manufacturing~\cite{Pe97}  or the design of envelope data structures
supporting ray shooting queries (seeking the first sphere, if any,
met by a query ray)~\cite{AAS97}.
In~\cite{MPT01,Me01,M2}, the question of arrangements of four spheres
in $\R^3$ with an infinite number of common tangent lines is discussed
from various viewpoints.
This was recently characterized~\cite{BGLP04}.

We first study the common tangents to two spheres that also meet a given line
and show that this 1-dimensional family (a curve) of lines determines the two
spheres.  
When the configuration of the two spheres and fixed line becomes degenerate, 
this curve becomes reducible and we classify its possible components.
In most cases, the components of the curve also determine the spheres and it
is only in the remaining highly degenerate cases that there is a third
sphere tangent to all the lines in a given component of the curve.

\begin{thm}\label{T:ls2}
 Let $\ell$ be a line and $S_1$, $S_2$ and $S_3$ be spheres in $\R^3$.
 Then there are infinitely many lines that meet $\ell$ and are 
 tangent to each sphere in precisely the following cases:
\begin{enumerate}
 \item [(i)] The spheres are tangent to each other at the same point 
             and either {\rm (a)} $\ell$ meets that point, 
             or {\rm (b)} it lies in the common tangent plane, or both. 
             The common tangent lines are the lines in the tangent plane
             meeting the point of tangency.
 \item [(ii)] The spheres are tangent to a cone whose apex lies on 
             $\ell$.
             The common tangent lines are the ruling of the cone.

 \item [(iii)] The spheres meet in a common circle and the line
             $\ell$ lies in the plane of that circle.
             The common tangents are the lines in that plane tangent to the
             circle.
 \item [(iv)] The centres of the spheres lie on a line 
             $m$ and $\ell$ is tangent to all three spheres.
             The common tangent lines are one ruling on the hyperboloid of
             revolution obtained by rotating $\ell$ about $m$. 
\end{enumerate}
 In cases {\rm (ib)}, {\rm (iii)} and {\rm (iv)} lines parallel to $\ell$ are
 excluded. 
\end{thm}

The classification in real projective space is similar, except that lines
parallel to $\ell$ are not excluded, and $\ell$ might be at infinity
in (ib) or (iii).
Also, in (ii), the cone may be a cylinder (cone with apex at infinity) with
$\ell$ parallel to the cylinder.
We do not know the classification over the complex numbers, as the real numbers
are used in an essential way to rule out several possibilities.
\begin{figure}[htb]
\[
  \begin{picture}(128,110)(-20,-5)
   \put(0,10){\epsfysize=90pt\epsffile{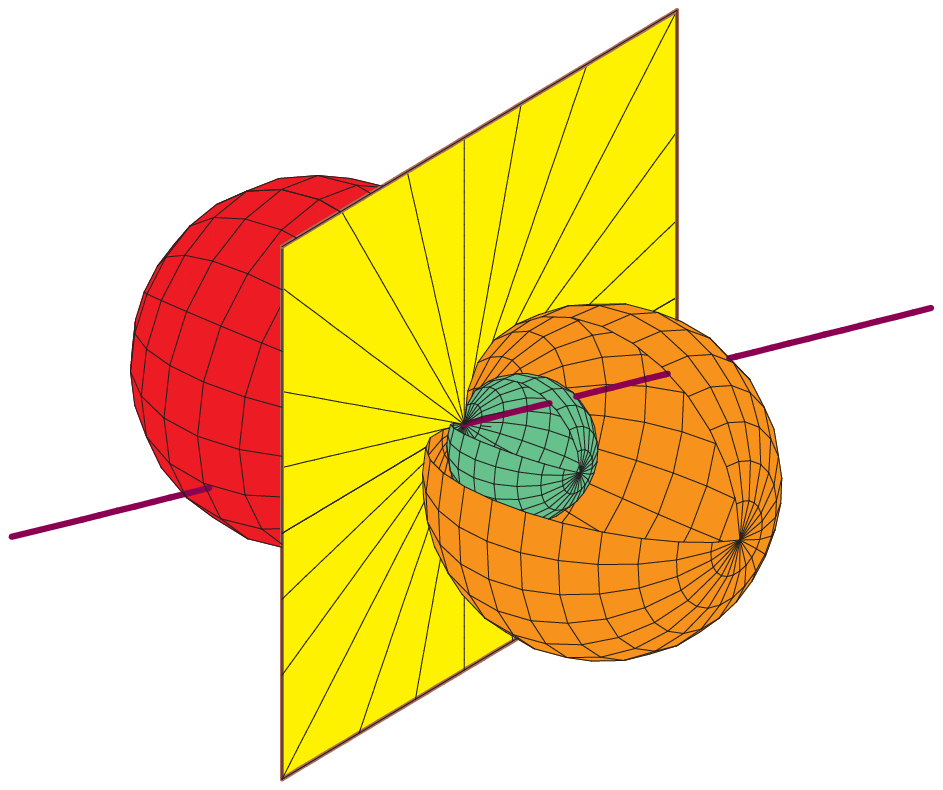}}
   \put(90,70){$\ell$}
   \put(40,-5){(ia)}
  \end{picture}
 \qquad
  \begin{picture}(85,110)(0,-5)
   \put(0,10){\epsfysize=90pt\epsffile{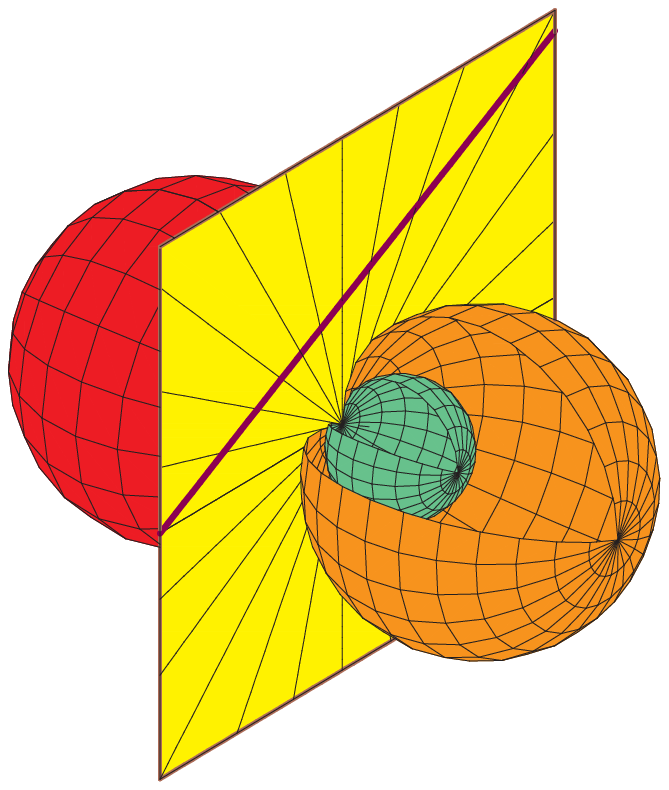}}
   \put(70,90){$\ell$}
   \put(35,-5){(ib)}
   \end{picture}
 \qquad
  \begin{picture}(175,110)(0,-5)
   \put(0,10){\epsfysize=90pt\epsffile{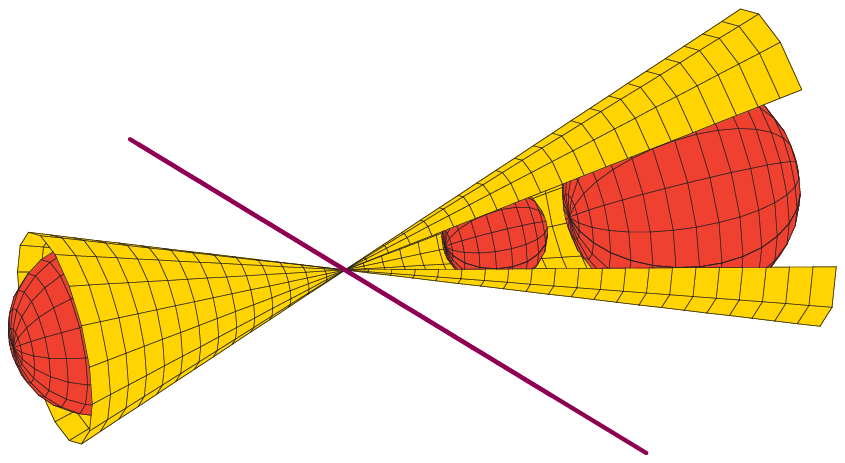}}
   \put(40,72){$\ell$}
   \put(85,-5){(ii)}
  \end{picture}
\]

\[
  \begin{picture}(105,110)(0,-10)
   \put(0,10){\epsfysize=100pt\epsffile{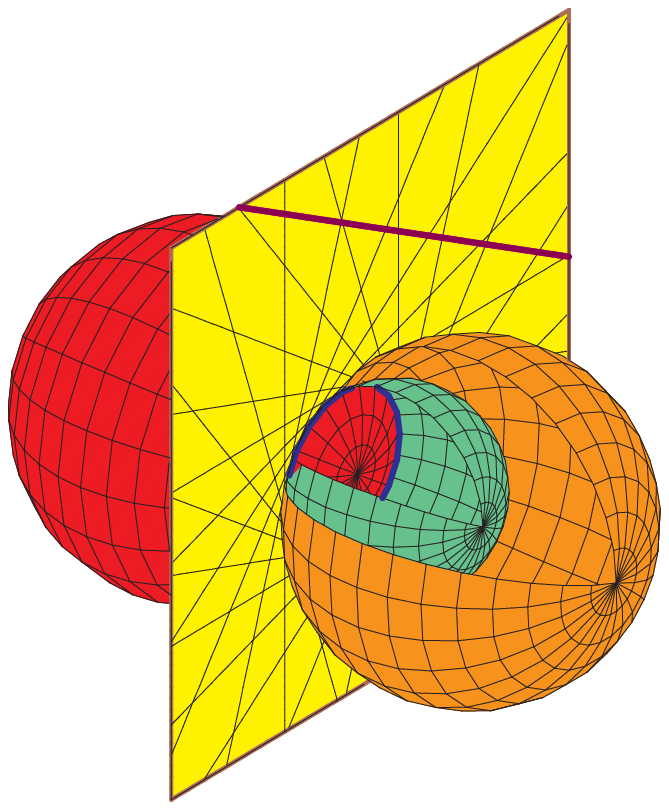}}
   \put(80,85){$\ell$}
   \put(50,-5){(iii)}
  \end{picture}
 \qquad
  \begin{picture}(249,83)(0,-10)
   \put(0,15){\epsfxsize=240pt\epsffile{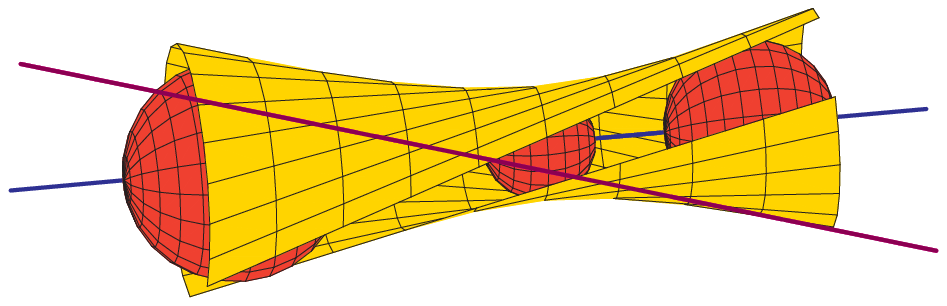}}
   \put(5,35){$m$}
   \put(15,82){$\ell$}
   \put(120,-5){(iv)}
  \end{picture}
\]    

\caption{The Geometry of Theorem~1.}\label{F:Thm1}

\end{figure}        

The primary difficulty in proving Theorem~\ref{T:ls2} lies not in recognising
the (fairly obvious) possibilities, but rather in excluding all others.
In Section 1 we give algebraic and geometric preliminaries concerning the
geometry of lines tangent to spheres.
In Section~\ref{Sec:Main} we study the envelope of lines tangent to two spheres
that meet a fixed line, describe the configurations of Theorem~\ref{T:ls2}
and begin to exclude other possibilities.
Specifically, we show that the algebraic subset $\tau$ of the Grassmannian
consisting of common tangents to two spheres and a line is a curve of degree at
most 8. 
We show that it is impossible for any component $\sigma$ of $\tau$ to have
degree 3 or 5 or 7, and if a component $\sigma$ has degree 4 or 6 or 8, then it 
determines the two spheres. 
Only in the cases of $\sigma$ having degree 1 or 2 can there be more than
2 spheres; these are described in Theorem~\ref{T:ls2}. 
The intricate argument when $\sigma$ has degree 4 is treated separately in
Section~\ref{S:Quartics}. 

\section{Real line geometry and Pl\"ucker coordinates}\label{S:Pluecker}

We review some aspects of the geometry of lines in space, Pl\"ucker
coordinates for lines, and preliminaries from real algebraic geometry.
Good general references are~\cite{CLO92,GH78,PW01}.

While our main concern is lines in $\R^3$ tangent to spheres, some algebraic
arguments we make will apply to lines in complex projective 3-space $\P^3$
tangent to a quadric surface, and this added generality will sometimes be
necessary.
Nevertheless, at key junctures the real-number nature of the answers we
seek will be used in an essential way.
For example, distinct concentric spheres have no common tangents in
$\R^3$, but do have a 2-dimensional family of common complex tangents in
$\P^3$. 
We use $\R^3$, $\R\mathbb{P}^3$ and $\mathbb{P}^3$ to indicate real affine space,
real projective space, and complex projective space, respectively.
In general, algebraic computation takes place in $\mathbb{P}^3$, and we then
interpret the answers in $\R^3$ to obtain our restrictions and to describe the
final configuration.

\subsection{Common tangents in a plane and through a point}
A useful warm-up is the following elementary
determination of the common tangents two spheres 
that lie in a plane and the common tangents through a point.

\begin{prop}\label{P:2spherePlane}
 Let $S_1$ and $S_2$ be spheres and\/ $\Pi$ be a plane $\mathbb{R}^3$ 
 such that $\Pi\cap S_1\neq \Pi\cap S_2$, as subsets of\/ $\mathbb{R}^3$.
 Then there are four (complex) common tangents to $S_1$ and $S_2$ that lie in
 $\Pi$, counted with multiplicity.
 Furthermore, if\/ $\Pi$ is not tangent to either sphere, then these common
 tangents determine the circles $\Pi\cap S_1$ and $\Pi\cap S_2$.
\end{prop}

\begin{proof}
 If $\Pi$ is not tangent to $S_1$, then $\Pi\cap S_1$ is a smooth conic whose
 tangents form a conic in the dual projective plane $\widehat{\Pi}$. 
 If $\Pi$ is tangent to $S_1$ at a point $q$ then the `tangents' (lines with a
 point of double contact) to $\Pi\cap S_1$ are the set of lines in $\Pi$
 through $q$, counted with multiplicity 2, a degenerate conic in
 $\widehat{\Pi}$.  
 Thus the common (complex) tangents to $\Pi\cap S_1$ and $\Pi\cap S_2$ are the
 intersection of two conics in $\widehat{\Pi}$, and so there are four, counted
 with multiplicity.  

 Suppose that $\Pi$ is not tangent to either sphere.
 Then the four common tangents determine the circles $C_1:=\Pi\cap S_1$
 and $C_2:=\Pi\cap S_2$:
 three distinct (complex) lines in $\Pi$ are common tangents to four circles.
 The fourth common tangent selects two of these four circles.
 This is clear when there are four distinct common tangents (see
 Figure~\ref{F:circle}). 
\begin{figure}[htb]
 \[
   \setlength{\unitlength}{.91666pt}
   \begin{picture}(150,94)
    \put(16,0){\epsfxsize=110pt\epsffile{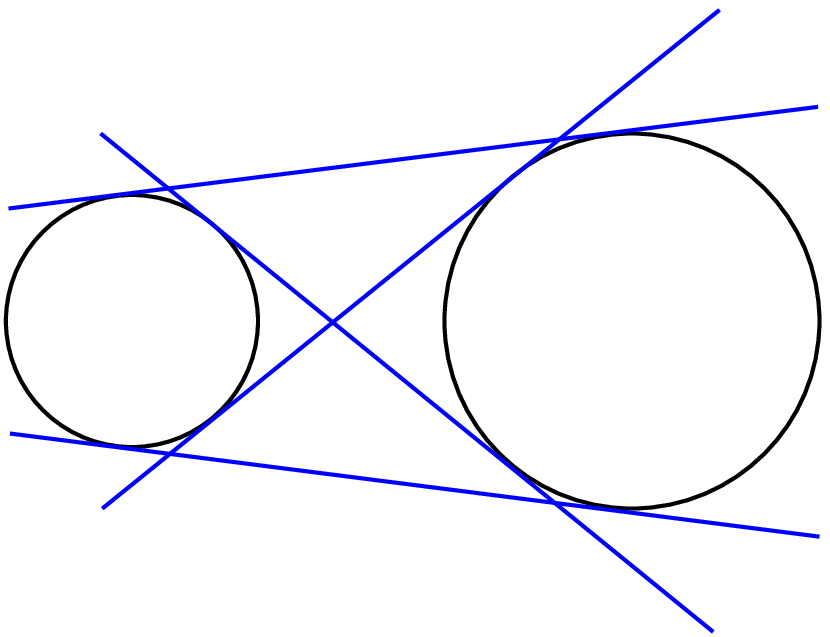}}
    \put(0,43){$C_1$}    \put(140,43){$C_2$}
   \end{picture}
   \qquad \qquad
   \begin{picture}(150,94)
    \put(17,0){\epsfxsize=105.4166pt\epsffile{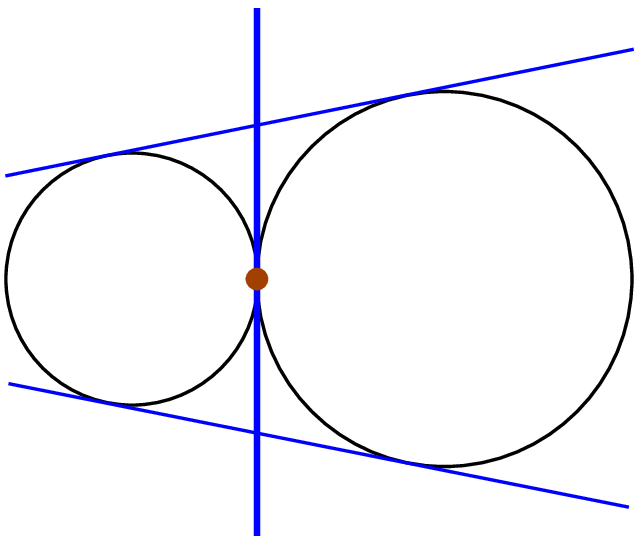}}
    \put(0,43){$C_1$} \put(138,43){$C_2$}
     \put(48,87){$m$} 
   \end{picture}
 \]
\caption{Lines tangent to two circles}
\label{F:circle}
\end{figure}
 Otherwise, $C_1$ and $C_2$ are tangent at the same point of the common
 tangent having multiplicity 2, and this additional information determines $C_1$
 and $C_2$. 
\end{proof}

\begin{prop}\label{P:cones}
 Let $S_1$ and $S_2$ be spheres and $p$ be a point in $\mathbb{R}^3$ such that
 the cone with apex $p$ tangent to $S_1$ is distinct from the cone with apex $p$
 tangent to $S_2$. 
 Then there are four (possibly complex) common tangents to $S_1$ and $S_2$ that
 meet $p$, counted with multiplicity.
\end{prop}

\begin{proof}
 Let $C_1$ and $C_2$ be the cones with apex $p$ tangent to $S_1$ and $S_2$,
 respectively. 
 If $\Pi$ is a plane not containing $p$,
 then $\Pi\cap C_1$ and $\Pi\cap C_2$ are plane conics, distinct by hypothesis. 
 These conics meet in four points, counted with multiplicity, and these points
 of intersection are where the common tangents meet $\Pi$.
\end{proof}

We display a configuration of two cones with four real common tangent lines below.
\[
  \epsfxsize=150pt\epsffile{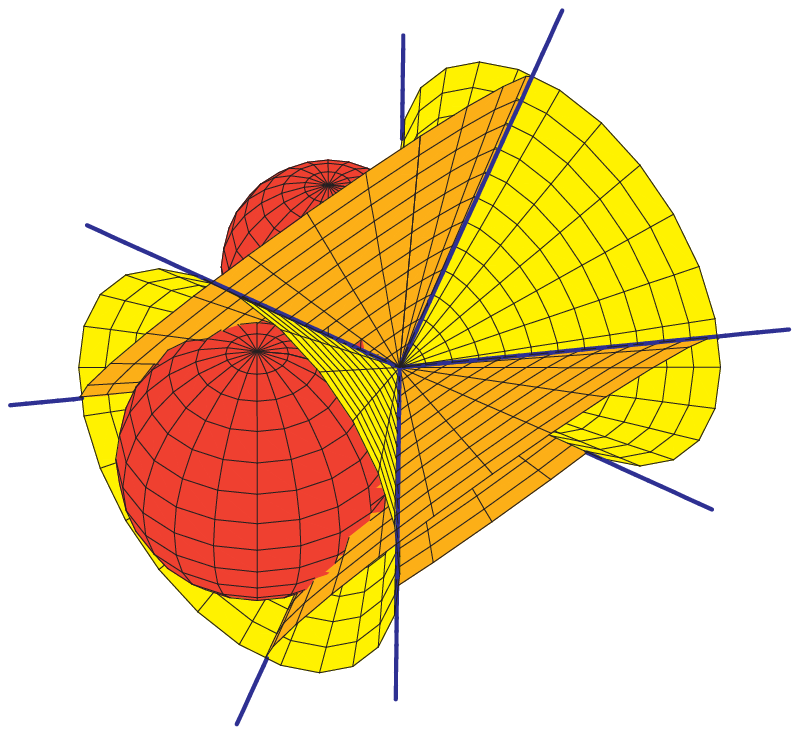}
\]

\subsection{Pl\"ucker coordinates for lines}
The Grassmannian $\G$ (or Klein quadric~\cite[\S 2.1]{PW01}) is the
set of lines in complex projective 3-space $\P^3$.
The line $\ell$ between two points $x=(x_0,x_1,x_2,x_3)^{\rm T}$ and 
$y=(y_0,y_1,y_2,y_3)^{\rm T}$ in $\P^3$ is represented (non-uniquely) by the  
$4\times 2$-matrix whose columns are the vectors $x$ and $y$.
A unique representation is given by its Pl\"ucker vector
$p_\ell:=(p_{01},p_{02},p_{03},p_{12},p_{13},p_{23})^{\rm T}\in\P^5$, where
$p_{ij}:=x_iy_j-x_jy_i$.
The coordinates of the Pl\"ucker vector satisfy the Pl\"ucker relation
 \begin{equation}\label{E:Pluecker}
  0\ =\ p_{03}p_{12}\ -\ p_{02}p_{13}\ +\ p_{01}p_{23}\,,
 \end{equation}
and $\G$ is the set of points in $\P^5$ which satisfy this relation.
Moreover, a line $\ell$ is in $\R\P^3$ if and only if its Pl\"ucker vector
lies in $\R\P^5\cap\G$, the real part of the Grassmannian.

A line $\ell$ meets another line $\ell'$ if and only if their Pl\"ucker
vectors $p$ and $p'$ satisfy
\begin{equation}
\label{E:linesMeet}
    p_{01} p'_{23}
  - p_{02} p'_{13} 
  + p_{03} p'_{12} 
  + p_{12} p'_{03} 
  - p_{13} p'_{02} 
  + p_{23} p'_{01}\ =\ 0\,.
\end{equation}
Geometrically, this means that the set of lines meeting a given line $\ell$
is the intersection of $\G$ with a hyperplane $\Lambda_\ell$ in $\mathbb{P}^5$.
For example, the $x$-axis (defined by the vanishing of the last 2 coordinates)
is the span of the points $(1,0,0,0)^{\rm T}$ and $(0,1,0,0)^{\rm T}$, and
thus has Pl\"ucker vector $(1,0,0,0,0,0)^{\rm T}$ as $p_{01}=1$ and the
other coordinates vanish.
The set of lines that meet the $x$-axis is then given by 
$   p_{23}=0.$

Each hyperplane $\Lambda_\ell$ is tangent to $\G$ at the point $p_\ell\in\G$.
The left-hand side of~\eqref{E:linesMeet} defines a symmetric bilinear form
$\langle p,p' \rangle$ on $\mathbb{P}^5$, identifying it with the dual
projective space of hyperplanes in $\mathbb{P}^5$, and $\G$ with its dual
$\G^*$. 

\begin{rem}\label{R:Three} 
 Let $H\subset\mathbb{P}^5$ be a 3-dimensional linear subspace.
 The set of hyperplanes containing $H$ is a line $h$ in the dual
 projective space.
 If the intersection $h\cap\G^*$ consists of two points, then these correspond
 to two hyperplanes $\Lambda_\ell$ and $\Lambda_{\ell'}$.
 Thus $H=\Lambda_\ell\cap\Lambda_{\ell'}$ and $H\cap\G$ parametrizes lines
 meeting $\ell$ and $\ell'$.
 If $H$ is real, then either both $\ell$ and $\ell'$ are real or else they are
 complex conjugate to each other. 
 If however $h\cap\G^*$ is a single point, so that
 $h$ is tangent to $\G^*$, 
 then $H$ is tangent to $\G$ and $H\cap\G$ is a singular quadric---a
 cone over a plane conic with apex $p_\ell$, the point of tangency. 
\end{rem}

Pl\"ucker coordinates afford a compact characterization of the lines
tangent to a given sphere and more generally tangent to a quadric surface.
We follow the presentation of~\cite{2l2s,STh02b}.
Identify a quadric $x^{\mathrm{T}} Q x = 0$ in $\mathbb{P}^3$ with the
symmetric  $4 \times 4$-matrix $Q$. 
Thus the sphere with centre $(x_0,y_0,z_0)^{\mathrm{T}}\in\mathbb{R}^3$ and
radius $r$ described in $\mathbb{P}^3$ by 
$(x - x_0 w)^2 + (y - y_0 w)^2 + (z - z_0 w)^2 = r^2 w^2$
is identified with the matrix
 \[
  \left(
    \begin{matrix}
      x_0^2 + y_0^2 + z_0^2 - r^2 & -x_0 & -y_0 & -z_0 \\
      -x_0 & 1 & 0 & 0 \\
      -y_0 & 0 & 1 & 0 \\
      -z_0 & 0 & 0 & 1 
    \end{matrix}
  \right) \, .
 \]

The second exterior power of a $4\times 4$-matrix $Q$ is the
$6\times 6$-matrix $\wedge^2Q$ whose entries are the $2\times 2$-minors of
$Q$.
Its rows and columns have the same index set as do Pl\"ucker vectors. 

\begin{prop} [Proposition 5.2 of~\cite{M2}]
\label{L:tangentCondition}
 A line $\ell \subset \mathbb{P}^3$ is tangent to a quadric $Q$ if and only
 if its Pl\"ucker vector $p_\ell$ lies on the quadric 
 hypersurface in $\mathbb{P}^5$ defined by $\wedge^2 Q$, 
 \begin{equation}
  \label{E:tangentCondition}
  p_\ell^{\mathrm{T}} \, \bigl(\wedge^2 Q\bigr) \, p_\ell\ =\ 0\,.
 \end{equation}
\end{prop}

For a sphere with radius $r$ and centre 
$(x_0, y_0, z_0)^T \in \mathbb{R}^3$ 
the quadratic form $p_\ell^{\mathrm{T}} \bigl(\wedge^2 Q\bigr) p_\ell$ is
 \begin{equation} \label{eq:tangentEquation}
  \left(
   \begin{array}{c}
    p_{01} \\ p_{02} \\ p_{03} \\ p_{12} \\ p_{13} \\ p_{23}
   \end{array}
   \right)^{\mathrm{T}}
   \left(
   \begin{array}{c@{~}c@{~}c@{~}c@{~}c@{~}c}
     y_0^2 + z_0^2 - r^2 & - x_0 y_0 & - x_0 z_0 & y_0 & z_0 & 0 \\
     - x_0 y_0 & x_0^2 + z_0^2 - r^2 & - y_0 z_0 & -x_0 & 0 & z_0 \\
     - x_0 z_0 & - y_0 z_0 & x_0^2 + y_0^2 - r^2 & 0 & -x_0 & -y_0 \\
     y_0 & -x_0 & 0 & 1 & 0 & 0 \\
     z_0 & 0 & -x_0 & 0 & 1 & 0 \\
     0 & z_0 & -y_0 & 0 & 0 & 1
   \end{array}
  \right)
  \left(
  \begin{array}{c}
    p_{01} \\ p_{02} \\ p_{03} \\ p_{12} \\ p_{13} \\ p_{23}
  \end{array}
  \right) \, .
 \end{equation}

We \rule{0pt}{14pt}use the relation between the intrinsic geometry of
$\G$ and families of lines in $\P^3$.
This begins with the following identification of linear subspaces lying in $\G$.

\begin{prop}[\cite{PW01}, pp.~141--143]\label{P:Pline}
 The set of lines in $\P^3$ meeting a given point $p$ is a $2$-dimensional
 linear subspace contained in $\G$.
 Dually, the set of lines lying in a given $2$-plane $\Pi$ is a $2$-dimensional
 linear subspace contained in $\G$, and any $2$-plane contained in
 $\G$ is one of these two types.
 
 A $1$-dimensional linear subspace contained in $\G$ consists exactly of the set
 of lines in $\P^3$ which contain a fixed point $p$ and lie in a fixed plane $\Pi$.
\end{prop}

\subsection{The envelope of a curve in $\G$}
\begin{figure}[htb]
\[
  \begin{picture}(190,160)
   \put(0,0){\epsfxsize=2.6in\epsffile{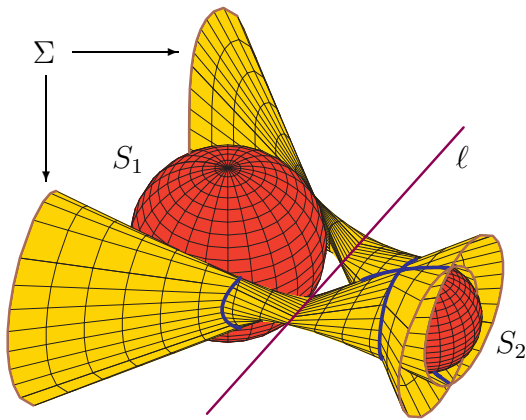}}
   \put(40,95){$S_1$}  \put(185,25){$S_2$}  \put(170,95){$\ell$}
   \put(10,135){$\Sigma$}\put(15,130){\vector(0,-1){40}}
                         \put(25,139){\vector(1,0){40}}
  \end{picture}
\]
\caption{An envelope}
\label{F:envelope}
\end{figure}

The union $\Sigma$ of the lines in $\mathbb{P}^3$ determined by a curve $\sigma$
in $\G$ is a ruled surface, which we call the {\it envelope} of $\sigma$.
The lines in $\sigma$ are a ruling of $\Sigma$.
For example, there is an irreducible octic (degree 8) curve $\sigma$ in $\G$
consisting of lines meeting the $x$-axis that are also tangent to two spheres
centred at $(0,\pm2,0)$ with respective radii $1$ and $\sqrt{3}$.
Its envelope $\Sigma$ is displayed in Figure~\ref{F:envelope}, where we 
also show the curves of self-intersection of $\Sigma$.
The real  points of the envelope $\Sigma$ of Figure~\ref{F:envelope} appear to
have two irreducible components, and thus contain two distinct families of
lines. 
This is a real-number illusion.
The complex curve is irreducible (in fact smooth and connected), but its real
points consist of two topological components which cannot be separated by an
algebraic function.

Suppose that $\sigma$ is an irreducible curve in $\G$.
Unless its envelope $\Sigma$ is a plane $\Pi$, 
the degree of the surface $\Sigma\subset\mathbb{P}^3$ equals the degree of
the curve $\sigma$~\cite[Theorem~5.2.8]{PW01}, and this common degree is
the intersection multiplicity of $\Sigma$ with a line $m$
not contained in $\Sigma$.
If a plane $\Pi$ is not contained in $\Sigma$, then the degree of $\Sigma$
equals the degree of the plane curve $\Sigma\cap\Pi$, if we keep track of
multiplicities of its components. 
These methods also determine the degree of a curve $\sigma$, even when
$\sigma$ is reducible or has components whose envelope is a plane, if we are
careful with the multiplicities. 

We make free use of continuity arguments from algebraic geometry.
Specifically,  the intersection multiplicity of a generic configuration is a
lower bound for the multiplicity of a special configuration.

\begin{ex}\label{E:conic}
 An important case of ruled surfaces is when $\sigma$ is an irreducible plane
 conic.
 If the plane of $\sigma$ lies in $\G$, then there are two possibilities, by
 Proposition~\ref{P:Pline}. 
 When the plane consists of lines through a point $p\in\mathbb{P}^3$, then
 $\Sigma$ will be a cone (over a conic) with apex $p$, or a cylinder if $p$ lies
 at infinity.
 When the plane consists of lines in a plane $\Pi\subset\mathbb{P}^3$,
 then $\sigma$ will be the set of lines tangent to a conic in $\Pi$.

 The most interesting case is when the plane of $\sigma$ does not lie in $\G$.
 A smooth quadric surface in $\mathbb{P}^3$ is ruled by two families
 of lines and each line in one family meets each line in the
 other, but no other lines in its own family.
 Any smooth conic $\sigma\subset\G$ whose plane does not lie in $\G$
 corresponds to a ruling of some quadric surface $Q$~\cite[Prop.~3.3.1]{PW01}.
 A plane tangent to a point of $Q$ meets $Q$ in two lines, one from each
 ruling.  

 Suppose that $Q\subset\mathbb{R}^3$ is a real quadric surface.
 If $Q$ contains even a single real line, then $Q$ is either a hyperbolic
 paraboloid or a hyperboloid of one sheet, and the two rulings are real.
 Otherwise, the rulings are complex conjugate to each other.
 Any curve $\gamma$ on $Q$ has a {\it bidegree} $(a,b)$, where $a$ is the
 intersection multiplicity of $\gamma$ with a general line in one family and $b$
 is the intersection multiplicity with a general line in the other.
 If $\gamma$ is the intersection of $Q$ with a hypersurface of degree $d$, then 
 its bidegree is $(d,d)$.
 Real curves $\gamma$ on a sphere $Q$ will have symmetric bidegree, as
 complex conjugation preserves $\gamma$ but interchanges the rulings of
 $Q$. 
\end{ex}

\subsection{Lines tangent to a sphere that meet a fixed line}\label{S:birat}
We study the correspondence between the points of a sphere and the lines tangent
to it that meet a fixed line.
Our goal is Lemma~\ref{L:degree}, which relates the degree of a curve $\gamma$
on the sphere to the curve $\sigma$ parametrizing lines which are 
tangent to the sphere at points of $\gamma$ and also meet the fixed line.

Let $S$ be a sphere and let $\ell$ be a line in $\mathbb{R}^3$.
The tangent plane at a general point $p$ of $S$ meets $\ell$ in a single point, and
thus there is a unique tangent line to $S$ at $p$ meeting $\ell$.
This defines an algebraic correspondence between points $p$ of $S$ and lines
tangent to $S$ that meet $\ell$.
This correspondence is given by polynomials that are quadratic in the
coordinates $(x_0,y_0,z_0)^T$ of the point $p$.

To see this, suppose that $S$ has centre $(a,b,c)^T\in\mathbb{R}^3$ and radius $r$.
Then the tangent plane to $S$ at a point $p=(x_0,y_0,z_0)^T$ is defined by the
linear equation
\[
   (x_0-a,\, y_0-b,\, z_0-c)\cdot (x,y,z)^T\ =\ 
    r^2 + a(x_0-a) + b(y_0-b) + c(z_0-c)\,.
\]
Write this compactly in vector form as $M\cdot(x,y,z)^T=\mu$,
where $M=(x_0{-}a, y_0{-}b,z_0{-}c)$ and 
$\mu= r^2 +M\cdot (a,b,c)^T$.

Suppose that the line $\ell$ is defined by the intersection of two linear
equations
\[
  \Lambda_1\cdot(x,y,z)^T\ =\ \lambda_1\qquad\textrm{and}\qquad
  \Lambda_2\cdot(x,y,z)^T\ =\ \lambda_2\,.
\]
Then the tangent plane to $S$ at $p$ meets $\ell$ at the solution of the linear
system
\[
   \left(\begin{array}{c}M\\ \Lambda_1\\ \Lambda_2\end{array}\right)
    \cdot
   \left(\begin{array}{c}x\\ y\\z\\\end{array}\right)
   \ =\ 
   \left(\begin{array}{c}\mu\\ \lambda_1\\ \lambda_2\end{array}\right)
   \ .
\]
This solution is 
 \begin{equation}\label{E:solution}
    \left(\begin{array}{c}x_1\\ y_1\\z_1\\\end{array}\right)
   \ =\ 
   \left(\begin{array}{c}M\\\Lambda_1\\\Lambda_2\end{array}\right)^{-1}
    \cdot
   \left(\begin{array}{c}\mu\\ \lambda_1\\ \lambda_2\end{array}\right)
   \ =\ 
   \textrm{Ad}\left(\begin{array}{c}M\\\Lambda_1\\\Lambda_2\end{array}\right)
    \cdot
   \left(\begin{array}{c}\mu\\ \lambda_1\\ \lambda_2\end{array}\right)
    \bigg/
   \det\left(\begin{array}{c}M\\ \Lambda_1\\ \Lambda_2\end{array}\right)\ .
 \end{equation}

The Pl\"ucker coordinates of the line in $\mathbb{P}^3$ tangent to $S$ at
$(x_0,y_0,z_0)^T$ that meets $\ell$ are the $2\times 2$ minors of 
\[
  \left[\begin{array}{cccc}1&x_0&y_0&z_0\\
                           1&x_1&y_1&z_1\end{array}\right]^T\ .
\]
By~\eqref{E:solution}, the coordinates $x_1,y_1,z_1$ are rational functions of
$x_0,y_0,z_0$ whose numerators are linear and which have the same (linear in
$x_0,y_0,z_0$) denominator, the determinant in~\eqref{E:solution}.
Clearing denominators shows that the Pl\"ucker coordinates of the line tangent
to $S$ at $(x_0,y_0,z_0)^T$ that meet $\ell$ are quadratic polynomials in
$x_0,y_0,z_0$. 

Let $\Phi$ be this algebraic correspondence between (complex) points of $S$ and
tangent lines 
 \begin{equation}\label{E:correspondence}
  \Phi\ \colon\ \{\textrm{points of $S$}\}\ \ - - \to\ \ 
                \{\textrm{lines tangent to $S$ meeting $\ell$}\}\,.
 \end{equation}
The broken arrow is used in algebraic geometry to indicate that $\Phi$ does not
define a function at all points of $S$.
Such a correspondence that is 1-1 almost everywhere is called a birational
equivalence.
We describe the loci where $\Phi$ is not 1-1 in geometric terms.
For this, we work in complex projective space.
This description is valid more generally for quadrics in $\mathbb{P}^3$.

A {\it special point} of $S$ is a point either lying on $\ell$ or whose tangent
plane contains $\ell$.
Every line tangent to $S$ at a special point meets $\ell$, and these are exactly
the points where $\Phi$ is not well-defined.
In particular, the correspondence $\Phi$ associates a special point to the
1-dimensional subspace of $\G$ consisting of those lines tangent to $S$ at the
special point. 
Conversely, if $m$ is a line lying in $S$ that meets $\ell$ (a 
{\it special line}), then $\Phi$ maps every point of $m$ to $m$.

If $\ell$ is tangent to $S$ at a point $p$, then $p$ is the only special point
and the two lines lying in $S$ in the tangent plane at $p$ are the only
special lines.
Otherwise $\ell$ meets $S$ in two points and there are two planes through $\ell$
tangent to $S$, so there are four special points $p_1,p_2,p_3,p_4$.
The four lines in $S$ that meet $\ell$ are the special lines.
Two special points are always real and the other two are a complex
conjugate pair, the special lines  always form two
complex conjugate pairs. 
In the terminology of algebraic geometry, we say that $\Phi$ blows up the
special points and blows down the special lines.

Let $\gamma\subset S$ be an irreducible curve, not a special line.
The proper transform $\Phi(\gamma)$ of $\gamma$ is the closure of the image of
$\gamma_0$ under $\Phi$, where $\gamma_0$ is $\gamma$ with any special points
removed. 
In Lemma~\ref{L:config}, we will need to know the degree of $\Phi(\gamma)$,
where $\gamma$ is rational.
This determination uses more machinery from algebraic
geometry than is needed by
other parts of this paper.

\begin{lemma}\label{L:degree} 
Let $\gamma\subset S$ be a rational curve, not one of the special lines.
Then the degree of the proper transform $\Phi(\gamma)$ of $\gamma$ is
\begin{enumerate}
 \item[(i)]
   ${\displaystyle  2\deg \gamma-\sum_{k=1}^4 {\rm mult}_{p_k}\gamma}$\ \ 
    if there are four special points $p_1, p_2, p_3, p_4$.
 \item[(ii)]
   ${\displaystyle  2\deg \gamma-2{\rm mult}_p \gamma}$, \ 
     if there is a unique special point $p$.\rule{0pt}{14pt}
\end{enumerate}
\end{lemma}

\begin{proof} As $\gamma$ is rational, there is a rational map 
$\nu\colon \P^1\to \gamma\subset S\subset\mathbb{P}^3$, which is an isomorphism
outside a finite set.  Let us choose homogeneous coordinates on $\P^1$.  Any
rational map $\phi\colon \P^1\to\P^n$ is defined by an $(n{+}1)$-tuple of
homogeneous polynomials of the same degree.  The polynomials
themselves are not uniquely defined, but their ratios are, and we can
multiply or divide them by a common factor. If we divide all of them
by their greatest common divisor to obtain an $(n{+}1)$-tuple of coprime
polynomials, then they do not all vanish at the 
same point of $\P^1$, so $\phi$ is defined everywhere, it is a
morphism. Furthermore, if $\phi$ is generically 1-1, i.e., it is 1-1
outside a finite set, then the degree of the coprime polynomials
defining $\phi$ is the degree of the image $\phi(\P^1)$ in $\P^n$.

This means that $\nu$ can be defined by coprime homogeneous polynomials of
degree $\deg \gamma$ on $\P^1$.  Let $\Phi_{ij}$, $0\le i< j\le 3$, be
the quadratic polynomials defining the map $\Phi$~\eqref{E:correspondence}.
Under $\nu$, they pull back to 
homogeneous polynomials $\Psi_{ij}$, $0\le i< j\le 3$, of degree
$2\deg \gamma$ on $\P^1$, which define the map $\Psi=\Phi\circ\nu$ whose
image agrees with $\Phi(\gamma)$.

The special points are exactly the points of $S$ where all the forms
$\Phi_{ij}$ ($0\le i< j\le 3$) vanish. 
If $\gamma$ does not pass through any of the special points, then the forms
$\Psi_{ij}$ ($0\le i< j\le 3$) do not all vanish at the same point of
$\P^1$, so they are coprime, and therefore the degree of 
$\Psi(\P^1)=\Phi(\gamma)$ is $2\deg \gamma$.
If $\gamma$ does pass through a special point, then all the
$\Psi_{ij}$ ($0\le i< j\le 3$)  vanish at the inverse image of that
point on $\P^1$, so they have a common factor which can be
cancelled. The degree of the image drops by the 
degree of the common factor, and it is this degree that we need to calculate.

The equation $\Phi_{ij}=0$ defines a curve on $S$ unless
$\Phi_{ij}$ is identically zero. Assume that this curve meets $\gamma$ at
a point $q$. Let $s\in\P^1$ be a point such that $\nu(s)=q$, and let
$\gamma'$ be the branch of $\gamma$ at $q$ corresponding to $s$. 
The intersection multiplicity of  $\Phi_{ij}=0$ with
$\gamma'$ is the same as the order of vanishing of $\Psi_{ij}$ at $s$.

In case (i), each polynomial $\Phi_{ij}$ ($0\le i< j\le 3$) vanishes at the four
special points. The multiplicity of a branch of $\gamma$ at a special
point $p_k$ is a lower bound for the intersection multiplicity of the curve
$\gamma$ at $p_k$, so each polynomial $\Psi_{ij}$  ($0\le i< j\le 3$) 
vanishes at least to that order at the corresponding point of
$\nu^{-1}(p_k)$. Therefore the polynomials $\Psi_{ij}$  ($0\le i< j\le 3$) have a
common factor of degree at least 
$\sum_{k=1}^4 {\rm mult}_{p_k}\gamma$. If we choose coordinates
such that $\ell$ is the $x$-axis and the centre of $S$ lies on the
$y$-axis, then one may check by direct calculation that the
curves defined by $\Phi_{02}=0$ and $\Phi_{13}=0$ are smooth at all of
the special points and have different tangent directions at each
special point. If $s$ is a point of $\P^1$ such that $\nu(s)$ is a
special point $p_k$ and $\gamma'$ is the corresponding branch of
$\gamma$ at $p_k$,
then the intersection multiplicity of $\gamma'$ at $p_k$ with at least one of
the curves $\Phi_{02}=0$ and $\Phi_{13}=0$ is exactly
${\rm mult}_{p_k}\gamma'$, so the greatest common divisor of the $\Psi_{ij}$ 
($0\le i< j\le 3$) vanishes at $s$ to order ${\rm mult}_{p_k}\gamma'$
exactly. Therefore the degree of  greatest common divisor of the $\Psi_{ij}$ 
($0\le i< j\le 3$) is $\sum_{k=1}^4 {\rm mult}_{p_k}\gamma$, which
proves the lemma for case (i). 

In case (ii) we can choose coordinates such that $\ell$ is the
$x$-axis\rule{0pt}{13pt} 
and $S$ has affine equation $x^2+(y-1)^2+z^2=1$. Then
$\Phi_{01}=(x^2-y^2-z^2)/2$ and $\Phi_{03}=x z$. They both define
curves with a simple node at the special point $p=(0,0,0)^T$, and their
tangent directions there are different. If $s\in\nu^{-1}(p)$ and
$\gamma'$ is the corresponding branch of $\gamma$ at $p$, then
the intersection multiplicity of $\gamma'$ at $p$ with at least one of the
curves $\Phi_{01}=0$ and $\Phi_{03}=0$ is exactly  
$2{\rm mult}_p \gamma'$, so at least one of  
$\Psi_{01}$ and $\Psi_{03}$ vanishes to
order exactly  $2{\rm mult}_p \gamma'$ at~$s$. Therefore the degree of
the greatest common divisor of the $\Psi_{ij}$ 
($0\le i< j\le 3$) is $2{\rm mult}_p \gamma$, which proves the lemma
for case (ii).
\end{proof}

The conclusion of Lemma~\ref{L:degree} also holds for a general curve $\gamma$,
but we only need the result for rational curves.
The proof for general curves is more sophisticated as we do not have a rational
parametrization. 
We may either work with a local parameter at points of the normalization of
$\gamma$, or consider the multiplicity of basepoints in the rational map $\Phi$.
\section{The envelope of lines meeting a fixed line and tangent to two
spheres} \label{Sec:Main}

Let $\ell$ be a line and let $S_1$ and $S_2$ be spheres in $\R^3$.
We assume that the spheres have infinitely many common real tangents, so that
in particular neither sphere contains the other.
We begin to describe the envelope of lines meeting $\ell$ that
are also tangent to two spheres, $S_1$ and $S_2$.
We also outline and begin our proof of
Theorem~\ref{T:ls2}.

Let $\tau$ be the algebraic subset of the Grassmannian $\mathbb{G}$ defined by
the linear equation $\Lambda_\ell$ and two quadratic equations of the 
form~\eqref{eq:tangentEquation} for lines to be tangent to $S_1$ and $S_2$.
Then $\tau$ parametrizes the  common tangents to the spheres $S_1$ and $S_2$
that also meet the fixed line $\ell$, retaining information about algebraic
multiplicities. 

\begin{thm}\label{T:tau} 
 $\tau$ is a curve of degree $8$ that determines the
 spheres $S_1$ and $S_2$.
\end{thm}

\begin{proof}
 Let $m$ be a line in $\mathbb{R}^3$ that meets $\ell$ in a point $p$.
 Then the common tangents to $S_1$ and $S_2$ that meet both $m$ and $\ell$
 are either the common tangents to $S_1$ and $S_2$ that lie in the plane
 $\Pi$ spanned by $\ell$ and $m$, or the common tangents to $S_1$ and $S_2$
 passing through the point $p$.
 When the line $m$ is chosen so the point $p$ and the plane $\Pi$ satisfy the
 hypotheses of Propositions~\ref{P:2spherePlane} and~\ref{P:cones} then there
 will be four of each type, counted with multiplicity.
 Thus the hyperplane $\Lambda_m$ in Pl\"ucker space meets $\tau$ in eight points,
 counted with multiplicity.
 This implies that $\tau$ is a curve of degree 8.

 By Proposition~\ref{P:2spherePlane}, given a plane $\Pi$ through
 $\ell$ that is not tangent to either sphere, the tangents to $S_1$ and $S_2$
 lying in $\Pi$ determine the circles $\Pi\cap S_1$ and $\Pi\cap S_2$.
 Letting the plane $\Pi$ vary shows that $\tau$ determines the spheres $S_1$
 and $S_2$.
\end{proof}

\begin{cor}\label{C:tau}
 There cannot be infinitely many tangent lines to three spheres that also meet a
 given line, if the line and any two of the spheres are in general position so
 that the set of common tangents $\tau$ to the spheres that meet $\ell$ is an
 irreducible curve.  
\end{cor}

By Corollary~\ref{C:tau}, if we have a third sphere $S_3$ tangent to infinitely
many of the lines in $\tau$, that set of lines is a proper subset of
$\tau$. 
Since this infinite set is an algebraic subvariety of the curve
$\tau$, we conclude that $\tau$ is reducible and the subset of $\tau$
consisting of lines tangent to $S_3$ is a union
of some, but not all, components of $\tau$, together with possibly finitely
many other points of $\tau$.
We study the envelope $T$ of the curve $\tau$, a surface of degree 8.
 
\begin{thm}\label{T:mult4}
 $T$ is singular along the line $\ell$.
 Furthermore, if\/ $\Pi$ is a plane through $\ell$, but not lying in $T$,
 then $\ell$ is a component of multiplicity at least $4$ in the degree $8$ curve
 $\Pi\cap T$.
 If $\ell$ is not tangent to both spheres, then this multiplicity is exactly $4$.
\end{thm}

\begin{proof}
 By Proposition~\ref{P:cones}, there will be four common complex tangents to the
 two spheres that meet a fixed point in $\mathbb{P}^3$.
 Thus when the line and two spheres are in general position, there are four lines
 in $T$ that meet at a general point $p$ of $\ell$.
 In particular, $\ell$ lies on $T$ and four 
 branches of $T$ meet along $\ell$, and so $T$ is singular along $\ell$.
 Since four branches of $T$ meet along $\ell$, if $\Pi$ is a plane through
 $\ell$ that is not contained in $T$, then $\ell$ is a component of multiplicity
 at least 4 in the degree 8 curve $T\cap\Pi$.
 The other components of $T\cap\Pi$ are common tangents to the circles
 $S_1\cap\Pi$ and $S_2\cap\Pi$. 
 Since there are fpur such common tangents counted with multiplicity, these common
 tangents contribute a multiplicity of at least 4 to $T\cap\Pi$.

 Thus, if $\ell$ is not tangent to both spheres, then both $\ell$ and the set of
 common tangents to the circles $S_1\cap\Pi$ and $S_2\cap\Pi$ each contributes a
 multiplicity of 4 to $T\cap\Pi$. 
%
%
\end{proof}

We classify the possible components $\sigma$ of $\tau$ that can occur.

\begin{thm}\label{T:sigma}
 Suppose that $\ell$, $S_1$ and $S_2$ are, respectively, a line and two
 spheres in $\R^3$, and that $\sigma$ is a reduced, irreducible 
 curve in the Grassmannian $\mathbb{G}$ such that every line in $\sigma$ meets  
 $\ell$ and is tangent to both of $S_1$ and $S_2$.
 Then $\sigma$ has degree either $1$, $2$, $4$, $6$ or $8$.
 Furthermore, if $\sigma$ has degree $1$, then it is a component 
 of multiplicity
 $2$ or $4$ in the curve $\tau$ of all 
 lines in the Grassmannian that meet $\ell$
 and are tangent to both $S_1$ and $S_2$.
\end{thm}

Our strategy to prove Theorem~\ref{T:ls2} is now clear:
We shall examine the geometry of the configurations of $\ell$, $S_1$ and
$S_2$ when the component $\sigma$ has degree 1, 2, 4 or 6, and in
each case describe completely the configurations of a third sphere tangent to
all the lines in $\sigma$.
The cases when $\sigma$ has degree 1 or 2 are detailed in Theorem~\ref{T:ls2},
and when the degree is $4$ or $6$, we show there are no possibilities for the
third sphere.
Theorem~\ref{T:tau} covers the case when the component has degree 8.

We prove Theorem~\ref{T:sigma} in this section.
The main point will be to exclude the possibilities of $\sigma$ having degree
3 or 5 or 7.
Along the way, we will also begin the proof of Theorem~\ref{T:ls2} by studying
in detail the geometry of the envelope of lines that meet  $\ell$
and are tangent to both $S_1$ and $S_2$.
As before, we first work in complex projective space to describe the possible
components algebraically, and then consider the situation in $\mathbb{R}^3$ to
further exclude components and to describe the configuration in $\mathbb{R}^3$.

\subsection{The degree of $\sigma$ is 1} 
We have the following trichotomy.

\begin{thm}\label{T:deg1}
 Let $\sigma$ be a $1$-dimensional linear subspace in $\G$ consisting of
 the lines lying in a plane $\Pi$ that meet a point $p\in\Pi$ such 
 that every line in $\sigma$ meets a fixed line $\ell$ and
 is tangent to two spheres $S_1$ and $S_2$.
 Then $S_1$ and $S_2$ are tangent to each other at the point $p$ with common
 tangent plane $\Pi$, and we further have one of
\[
   ({\rm i})\ p\in\ell\not\subset\Pi, \qquad 
   ({\rm ii})\ p\not\in\ell\subset\Pi \qquad or \qquad 
   ({\rm iii})\ p\in\ell\subset\Pi\,.
\]
 In {\rm (iii)}, $\sigma$ is a component of multiplicity $4$ in the degree $8$
 curve $\tau$ of common tangents to $S_1$ and $S_2$ that meet $\ell$, and in
 the remaining cases $\sigma$ is a component of multiplicity
 $2$.
\end{thm}

Theorem~\ref{T:deg1} proves statement (i) of Theorem~\ref{T:ls2}:
Suppose that $\sigma$ is a 1-dimensional linear subspace in $\G$, every line of
which meets $\ell$ and is tangent to three spheres $S_1$, $S_2$ and $S_3$.
By Theorem~\ref{T:deg1}, the three spheres are mutually tangent at the point
$p$ with common tangent plane $\Pi$, where $\sigma$ consists of the lines in
$\Pi$ that meet $p$, and we further have that $\ell$ either meets $p$ or $\ell$
lies in $\Pi$, or both.

\begin{proof}
 Our proof of this statement is entirely algebraic, working in complex
 projective space and in the Grassmannian $\G$.
 Suppose that $\sigma$ is a 1-dimensional linear subspace lying in $\G$.
 By Proposition~\ref{P:Pline}, $\sigma$ consists of the lines in a plane $\Pi$
 that contain a given point $p\in\Pi$.
 If each line in $\sigma$ is tangent to a sphere $S$, then each line is tangent
 to the conic $S\cap\Pi$.
 The only possibility is that $S$ is tangent to $\Pi$ at the point $p$.
 Similarly, if every line in $\sigma$ meets a fixed line $\ell$, then either
 (i) $\ell$ contains the point $p$ or else (ii) it lies in the plane $\Pi$ or
 (iii) both.

 We determine the multiplicity of $\sigma$ in each of these three cases.
 Suppose that $p$ is the origin in $\mathbb{R}^3$ and $\Pi$ is the
 $yz$-plane, so that $\sigma$ is defined by the equations 
 $p_{01}=p_{12}=p_{13}=p_{23}=0$.
 Then $S_1$ and $S_2$ have their centres along the $x$-axis, say at points $x_1$
 and $x_2$ (with radii $|x_1|$ and $|x_2|$).
 The quadratic forms~\eqref{eq:tangentEquation} defined by $\wedge^2S_i$ become
\[
   -x_i^2p_{01}^2 - 2x_i(p_{02}p_{12}+p_{03}p_{13})
   + p_{12}^2 + p_{13}^2 + p_{23}^2\,, \quad\mbox{for}\ i=1,2\,.
\]
 We see that the quadric hypersurfaces defined by $\wedge^2S_1$ and
 $\wedge^2S_2$ have the same tangent plane at 
 every point $(0,p_{02},p_{03},0,0,0)^T$ of $\sigma$ and therefore $\sigma$ has
 multiplicity at least 2 in $\tau$.
 We show that this multiplicity is exactly 2 in cases (i) and (ii) by considering
 degenerate configurations. 

 In case (i) ($p\in\ell$) suppose that $\ell$ is the $x$-axis and let $\tau$ be
 the curve of common tangents to $S_1$ and $S_2$ that also meet $\ell$.
 We compute the intersection of $\tau$ with the hyperplane $\Lambda_m$ in
 Pl\"ucker space corresponding to a line $m$ meeting $\ell$ in a point $q$
 distinct from the origin $p$ and from the apex of the cone tangent to both
 $S_1$ and $S_2$.
 A line meets both $\ell$ and $m$ if it lies in the plane $\Pi$ they span or if
 it meets $q$. 
 As in the proof of Theorem~\ref{T:tau}, there will be four common tangents to
 $S_1$ and $S_2$ meeting $q$ (but none lies in $\sigma$) and four common
 tangents lying in $\Pi$. 
 Since $\Pi\cap S_1$ and $\Pi\cap S_2$ are tangent at $p$, the common tangent
 at $p$ has multiplicity 2.
 Since only this common tangent is a line in $\sigma$, the 
 multiplicity of $\sigma$ in $\tau$ is at most 2, and hence equal to 2 in case
 (i). 

 In case (ii) ($\ell\subset\Pi$) we similarly suppose that $\ell$ is any line
 in the $yz$-plane not containing the origin $p$ and that $m$ is a line not
 lying in the $yz$-plane that meets $\ell$ at a point $q$.
 None of the four common tangents to $S_1$ and $S_2$ in the plane spanned by
 $\ell$ and $m$ is a line in $\sigma$, and there are four common tangents
 (counted with multiplicity) containing $q$.
 As before, there will be three such common tangents with the line
 $\overline{pq}$ having multiplicity 2.
 This completes the proof in case (ii) as $\overline{pq}$ is in $\sigma$.

 The proof in case (iii) is a consequence of Example~\ref{E:mult4} below.
\end{proof}

\begin{ex}\label{E:mult4}
 Suppose that the spheres $S_1$ and $S_2$ are tangent at the point
 $p=(0,0,0)^T$ with common tangent plane the $yz$-plane.
 Scaling the coordinates, we may assume that $S_1$ has centre $(1,0,0)^T$ with
 radius 1 and $S_2$ has centre $(-r,0,0)^T$ with radius $|r|$.
 We forbid the values of $-1,0,1$ for $r$:
 If $r=-1$, then the spheres coincide, $r=0$ is a degenerate sphere, and if
 $r=1$, then there is additional symmetry which we discuss in Section~\ref{S:4.2}.
 We also assume that $\ell$ is the $z$-axis.
 Then the curve $\tau$ of lines tangent to $S_1$ and $S_2$ that also
 meet $\ell$ is cut out by the polynomials $p_{12}$, and 
\begin{equation}\label{eq:mult4}
 \begin{array}{c}
  p_{03}p_{12}\ -\ p_{02}p_{13}\ +\ p_{01}p_{23}\,,\\
  -   p_{01}^2 - 2(p_{02}p_{12}+p_{03}p_{13})
                           +p_{12}^2+p_{13}^2+p_{23}^2\rule{0pt}{14pt}\,,
    \makebox[0pt][l]{\quad and} \\
  -r^2p_{01}^2 + 2r(p_{02}p_{12}+p_{03}p_{13})
                           +p_{12}^2+p_{13}^2+p_{23}^2\rule{0pt}{14pt}\,.
 \end{array}
\end{equation}

 Let $\sigma\subset\G$ be the set of
 lines in the $yz$-plane through the origin.
 As before, $\sigma$ is a component of the curve $\tau$ of common
 tangents to the spheres that meet the $z$-axis.
 The remaining components have degree at most seven, when counted with their
 multiplicities. 
 We invite the reader to check that under the following substitution
 \begin{equation}\label{E:parametrization}
   \begin{array}{ll}
      p_{01}\ =\ 2\sqrt{r}(s^2 + t^2)(s^2 - t^2)\ \ \ \ &
      p_{12}\ =\ 0\\
      p_{02}\ =\  4\sqrt{r}st(s^2 + t^2)&\rule{0pt}{14pt}
      p_{13}\ =\ 2r(s^2 - t^2)^2\\
      p_{03}\ =\ (r{-}1)(s^2 + t^2)^2\quad\rule{0pt}{14pt}& 
      p_{23}\ =\ 4rst(s^2 - t^2)
    \end{array}
\end{equation}
 the polynomials~\eqref{eq:mult4} vanish identically for all
 $[s,t]\in\mathbb{P}^1$. 
 The polynomials in $s,t$ in~\eqref{E:parametrization} have no common factor, so 
 this substitution defines a parametrized rational quartic curve $\rho$ in
 $\G$ which is a component of $\tau$ having degree 4.
 The bound of 7 for the degrees of the components of $\tau-\sigma$ shows that
 $\rho$ has multiplicity 1 in $\tau$.
 Since the five quartic polynomials in $s,t$ of~\eqref{E:parametrization} are
 linearly independent, $\rho$ is a rational normal quartic curve in the
 hyperplane $p_{12}=0$ of Pl\"ucker space.

 We claim that $\tau=\sigma\cup\rho$, which implies that $\sigma$ has
 multiplicity 4 in $\tau$.
 To prove this, we first observe that the lines in $\rho$ are the
 column space of the matrix 
 \begin{equation}\label{E:matrix}
  \left[
   \begin{array}{cccr}
     s^2+t^2 & 0 &0& -\sqrt{r}(s^2-t^2)\\
        0    & 2\sqrt{r}(s^2-t^2)& 4\sqrt{r}st &  (r{-}1)(s^2+t^2)
       \rule{0pt}{14pt}
  \end{array}
  \right]^T\ .
 \end{equation}
 Note that no line of $\rho$ coincides with the $z$-axis.

 Let $T$ be the envelope of $\tau$.
 Suppose that $\Pi$ is any plane containing the $z$-axis, other than the
 $yz$-plane which is the common tangent plane to $S_1$ and $S_2$.
 Then the intersection $T\cap\Pi$ consists of the $z$-axis and two
 remaining common tangents to the circles $\Pi\cap S_1$ and $\Pi\cap S_2$.
 We may suppose that $\Pi$ is defined by an equation $x=Ay$.
 By~\eqref{E:matrix} the lines in $\rho$ lying in $\Pi$ have parameter values
 satisfying $2\sqrt{r}(s^2-t^2)=4A\sqrt{r}st$, so there are exactly two lines in
 $\rho$ that lie in $\Pi$.
 Furthermore, none of these lines coincides with the $z$-axis.
 Since any line in $\tau$ that does not lie in $\sigma$ will span such a plane
 $\Pi$ with the $z$-axis, we see that lines in $\tau$ not in $\sigma$ must lie
 in $\rho$, which completes the proof that $\sigma$ has multiplicity 4 in
 $\tau$.

 In Figure~\ref{F:RatNorm}, we display the envelope $R$ of the
 rational normal curve $\rho$ when $r=2$.
 The envelope has been cut away to reveal both spheres.
 Drawn on the envelope are the two curves where the envelope is tangent to the
 spheres.
 The remaining arc of a circle is one component of the locus of
 self-intersection of the envelope, the other being the $z$-axis.
 \begin{figure}[htb]  
  \[
    \epsfxsize=3in\epsfbox{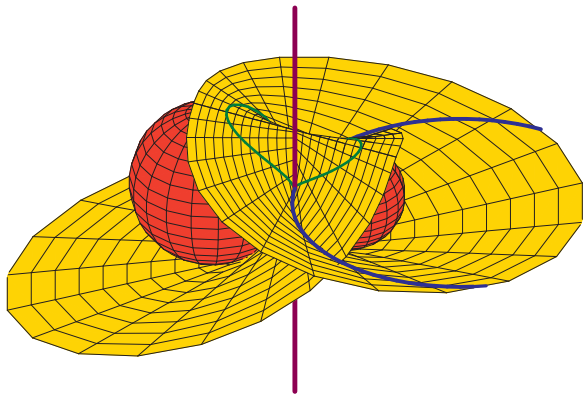}
  \]
 \caption{An envelope}\label{F:RatNorm}
\end{figure}
\end{ex}

\subsection{The degree of $\sigma$ is 2}\label{S:conic}
Suppose now that the curve $\tau$ of common tangents to $S_1$ and $S_2$ that
meet $\ell$ has a component $\sigma$ of degree 2.
As in Example~\ref{E:conic}, there are three possibilities for the envelope
$\Sigma$ of $\sigma$.
We discuss possible configurations of the spheres and line for each; these
possibilities correspond to cases (ii), (iii) and (iv) of
Theorem~\ref{T:ls2}.

\begin{enumerate}

\item[(i)] The envelope $\Sigma$ is a cone or cylinder.
     Because the lines on $\Sigma$ are tangent to the spheres $S_1$ and
     $S_2$, $\Sigma$ is circular and its axis $m$ contains the centres of the
     spheres.
     Since no line in $\mathbb{R}^3$ meets more than two lines on a circular
     cylinder, $\Sigma$ must be a cone with apex $p$ and $p\in\ell$.
     This gives case (ii) of Theorem~\ref{T:ls2}: $S_3$ is any sphere inscribed
     in the cone $\Sigma$ whose centre (on $m$) is any point distinct from the
     apex or the centres of $S_1$ and $S_2$. 

\item[(ii)] The envelope $\Sigma$ is the set of lines in a plane $\Pi$ tangent to
     a smooth conic $C$ in $\Pi$.
     Necessarily $\ell\subset\Pi$ and $\Pi\cap S_1=\Pi\cap S_2=C$, and $C$ is a
     circle.
     This gives case (iii) of Theorem~\ref{T:ls2}: $S_3$ is any other sphere with
     $\Pi\cap S_3=C$.

\item[(iii)] The envelope $\Sigma$ is a smooth quadric surface $Q$ with $\sigma$
     one of the rulings of $Q$.
     Then $\ell$ is a line in the other ruling.
     Since $Q$ contains the real line $\ell$, it is either a hyperbolic
     paraboloid or a hyperboloid of one sheet.
     As the lines in one ruling of $Q$ are tangent to the spheres $S_1$ and
     $S_2$, it cannot be a hyperbolic paraboloid.
     By Lemma~17 of~\cite{2l2s}, $Q$ is tangent to $S_1$ and $S_2$
     and is a hyperboloid of revolution with axis $m$ containing the centres of
     $S_1$ and $S_2$.
     Then $\ell$ is tangent to both $S_1$ and $S_2$ and $Q$ is obtained by
     revolving $\ell$ around the line $m$ spanning the centres of $S_1$ and
     $S_2$. 
     This gives case (iv) of Theorem~\ref{T:ls2}:  $S_3$ is any sphere
     inscribed in $Q$ whose centre (on $m$) is any point different than the
     centres of $S_1$ and $S_2$. 
\end{enumerate}

The following proposition will be useful in Sections~\ref{S:six} and
\ref{S:Quartics}. 

\begin{prop}\label{T:deg2}
 A component of $\tau$ of degree $2$ has multiplicity $1$ in $\tau$.
 There is a unique component of degree $2$ or of degree $1$ with multiplicity 
 $2$ unless either the centres of $S_1$ and $S_2$ both lie on $\ell$ or $S_1$
 and $S_2$ are symmetric with respect to $\ell$. 

 If the centres of $S_1$ and $S_2$ both lie on $\ell$, then $\tau$ has two
 components of degree $2$, one of which degenerates into a component of degree
 $1$ with multiplicity $2$ if $S_1$ and $S_2$ are tangent to each other.
 If we choose coordinates so that $\ell$ is the $x$-axis, then $\tau$ also
 contains the $1$-dimensional linear subspaces 
 \begin{equation}\label{E:special_lines}
   p_{01}\ =\ p_{23}\ =\ p_{02}\mp ip_{03}\ =\ p_{12}\mp ip_{13}\ =\ 0\,,
 \end{equation}
 each with multiplicity $2$.
 These parametrize lines in $\P^3$ that pass through a point of
 $\ell$ and one of the points at infinity $(0,0,1,\pm i)^T$ where the spheres
 are tangent to each other.
\end{prop}

Two spheres in $\mathbb{P}^3$ meet in two conics.
One is the imaginary spherical conic at infinity defined by $x_0=0$ and
$x_1^2+x_2^2+x_3^3=0$, and the other lies in a plane perpendicular
to the line segment joining the centres of the spheres.\medskip

\noindent{\it Proof.}
 {}From the geometric description of the configuration of the spheres and
 of $\ell$ when there is a component of multiplicity 2, there can
 only be one component of degree 2 or of degree 1 with multiplicity 2 in $\tau$,
 unless we are in one of two special cases:
 either (1) the centres of both spheres lie on $\ell$, or 
 (2) $\ell$ lies in a plane (in $\mathbb{P}^3$) in which the spheres meet and
 $\ell$ also meets the apex of a cone tangent to both spheres.
 Since $\ell\subset\R^3$, this plane is perpendicular to the line
 segment joining the centres of the spheres and moreover,  
 the spheres are necessarily symmetric about $\ell$. 
 In this case exactly one of
 the two components of degree 2 is real as either they meet in an imaginary
 conic and the  cone tangent to the spheres with apex in the plane of symmetry
 is real, or vice versa. 
 Also in this latter case, the remaining common tangents that meet $\ell$ form
 an irreducible quartic, illustrated in Figure~\ref{F:symmetric} and described
 in Theorem~\ref{T:quartics}.

 Consider the first special case.
 Suppose that both centres lie on $\ell$, which we take to be the $x$-axis, and
 consider the intersection of $S_1$ and $S_2$ with a real plane $\Pi$ through
 $\ell$.  
 These circles will have
 four common tangent lines consisting of two pairs symmetric with
 respect to $\ell$.  
 Rotating these tangent lines about $\ell$ gives two cones or cylinders, one of
 which degenerates into a plane if $S_1$ and $S_2$ are tangent to each other. 
 This gives two degree 2 components of $\tau$, or a degree 2 component and a 
 degree 1 component with multiplicity 2. 
 Such spheres are tangent to each other at the point $(0,0,1,\pm i)^T$ at
 infinity, with tangent plane containing $\ell$: this explains the two degree 1
 components of $\tau$ having multiplicity~2. 
 
 Now we show that a component $\sigma$ of degree 2 always has
 multiplicity 1. 
 By Proposition~\ref{P:Pline}, the lines lying in a plane $\Pi$ through 
 $\ell$ or passing through a point $p\in\ell$ form a 2-plane in $\G$. 
 In general, by Propositions~\ref{P:2spherePlane} and~\ref{P:cones}, there are
 four common tangent lines to $S_1$ and $S_2$ lying in $\Pi$ or passing 
 through any point $p$ of $\ell$, so the corresponding 2-planes in $\G$ meet
 $\tau$ at four points counted with multiplicity.

 Consider the three possibilities for a degree 2 component $\sigma$ of $\tau$.
 \begin{enumerate}

  \item[(i)] When the envelope $\Sigma$ of $\sigma$ is a cone, two lines of
    $\sigma$ lie in any plane $\Pi$ through $\ell$, but there will be two
    other common tangents.
    This accounts for the four lines of $\tau$ lying in the 2-plane
    $\Pi$. 
    Thus $\sigma$ must have multiplicity 1, for otherwise both lines of $\sigma$
    would have multiplicity 2 and there would be no others.

  \item[(ii)] When the envelope consists of tangents to a conic 
    $\Pi\cap S_1=\Pi\cap S_2$ in a plane $\Pi$ through $\ell$, two tangents in 
    $\sigma$ meet a general point $p$ of $\ell$, but there will be  
    two other common tangents to $S_1$ and $S_2$ meeting $p$.
    As before, this implies that $\sigma$ has multiplicity 1.

  \item[(iii)] Now suppose that $\ell$ is tangent to the spheres and $\sigma$ is
    one ruling of a hyperboloid of revolution.
    Every line $m$ of $\sigma$ meets $\ell$ and with it spans a plane $\Pi$.
    Unless $m\cap\ell$ is a point on one of the spheres, then $\Pi$ will contain
    four distinct common tangents to $S_1$ and $S_2$, thus two in addition to $m$ and
    $\ell$.
    Again, this implies that $\sigma$ has multiplicity 1.
\hfill $\Box$

 \end{enumerate}

\subsection{The degree of $\sigma$ is 3}

We show that this possibility does not occur.

\begin{prop}\label{P:cubic}
 The curve $\tau$ of common tangents to two spheres that meet a fixed line
 cannot have an irreducible component of degree $3$.
\end{prop}

\begin{proof}
 Let $\sigma\subset\tau$ be a component of degree 3.
 Then either $\sigma$ is a plane cubic or else it spans a 3-dimensional linear
 subspace of $\mathbb{P}^5$.
 Suppose that $\sigma$ is a plane cubic. 
 As $\G$ is a quadric, it can only contain $\sigma$ if it contains the whole
 plane it spans.
 This plane corresponds either to all lines through a given point in $\P^3$ or
 to all lines contained in a given plane, but in both cases the lines tangent to
 a sphere in such a plane are parametrized by a conic, not a cubic,
 and so this case cannot occur. 

 Now suppose that $\sigma$ spans a 3-dimensional subspace $H$ of $\mathbb{P}^5$.
 Note that $H\subset\Lambda_\ell$.
 There are two possibilities by Remark~\ref{R:Three}: either $H\cap\G$ is the
 set of lines meeting $\ell$ and another (uniquely defined) line $\ell'$, or
 $H\cap\G$ is a cone with apex $p_\ell$ over a plane conic.
 In the first case, both $\ell$ and $\ell'$ are real and $\sigma$ is a cubic
 curve of lines meeting $\ell$ and $\ell'$ that are also tangent to
 two spheres. 
 As shown in Section~5.1 of~\cite{2l2s}, this is impossible; the lines in a
 cubic in $\Lambda_\ell\cap\Lambda_{\ell'}$ are tangent to at most one sphere. 
 That determination was algebraic, and did not appeal to the real numbers.

 Consider the last possibility, that $H\cap\G$ is a cone with apex $p_\ell$
 over a plane conic. 
 Since $\sigma$ is a cubic curve in $H\cap\G$, is must meet the apex and so
 $\ell$ is tangent to both spheres.
 Let $\Sigma$ be the envelope of $\sigma$, a cubic surface in $\mathbb{P}^3$.
 The curve $\gamma\subset S_1$ along which lines in $\sigma$ are tangent to
 $S_1$ is a component with multiplicity 2 of the $(3,3)$-curve $\Sigma\cap S_1$.
 This implies that $\gamma$ is a $(1,1)$-curve, and that $\Sigma\cap S_1$
 contains components other than $\gamma$.
 These other components necessarily are the lines in $\sigma$ that lie in
 $S_1$. 
 Such a line must contain the point where $\ell$ is tangent to $S_1$, and it is
 not real, so $\sigma$ contains its complex conjugate.
 As both meet $\ell$ and $\ell\not\subset S_1$, they meet $\ell$ at its point of
 tangency with $S_1$.  
 Furthermore these lines and $\ell$ lie in the tangent plane to $S_1$ at that
 point. 
 This implies that $S_2$ is also tangent to $S_1$ at this same point.
 This configuration was studied in Example~\ref{E:mult4}, where it was shown
 that the irreducible components of $\tau$ have degrees 1 and 4, and not 3.
 This completes the proof of the impossibility of a cubic component.
\end{proof}

\subsection{The degree of $\sigma$ is 4}
In Section~\ref{S:Quartics} we show there are at most two spheres tangent to all
the lines in a degree 4 curve of lines meeting a fixed line.

\subsection{The degree of $\sigma$ is 5}

This possibility also does not occur:
if a component $\sigma$ of the curve $\tau$ of common tangents to the spheres
$S_1$ and $S_2$ that also meet $\ell$ has degree 5, then the residual curve
consisting of the other components of $\tau$ has degree 3.
By Proposition~\ref{P:cubic}, this degree 3 curve cannot be
irreducible and reduced. 
Theorem~\ref{T:deg1} excludes all the
possibilities that the degree 3 curve could split into lower degree
components. 

\subsection{The degree of $\sigma$ is 6}\label{S:six}

We show that in most cases a union of components of $\tau$ having degree 6
determines the spheres $S_1$ and $S_2$.

\begin{thm}\label{T:deg6} 
 Let $\tau$ be the curve of common tangents to spheres $S_1$ and $S_2$ that also
 meet a fixed line $\ell$ and suppose that $\tau$ is reducible with a
 component $\rho$ either of degree $2$, or of degree $1$ with multiplicity $2$. 
 Then the residual sextic $\sigma$ determines $S_1$ and $S_2$, except 
 when their centres both lie on $\ell$.
\end{thm}

By Proposition~\ref{T:deg2}, if the centres of $S_1$ and $S_2$ lie on $\ell$,
then $\tau$ has four components so that the residual sextic is reducible.

\begin{proof}
 Let $R$, $\Sigma$ and $T$ be, respectively, the envelopes of the curves
 $\rho$, $\sigma$ and $\tau$.
 Except when the centres of $S_1$ and $S_2$ lie on $\ell$ (which may be
 determined from $\Sigma$), we argue that we can determine the spheres $S_1$ and
 $S_2$ from $\sigma$.

 In Section~\ref{S:conic}, we discussed three possibilities for a component of
 $\rho$ with degree 2.
 The two cases (i) and (ii) of Theorem~\ref{T:deg1} for a component with degree
 1 and multiplicity 2 are, respectively, degenerations of cases (i) and (ii)
 of Section~\ref{S:conic}.
 There will be five cases to consider for the geometry of $\rho$ and we
 treat each case separately.

 In the first paragraph of each case we describe the
 number of lines in $\sigma$ meeting a general point of $\ell$ and the number
 of lines in $\sigma$ contained in a general plane $\Pi$ through $\ell$,
 together with some additional special geometry of that configuration.
 These descriptions will be different, showing that the sextic $\sigma$
 determines each case.
 Recall from Propositions~\ref{P:2spherePlane} and~\ref{P:cones} that $\tau$
 contains four lines through a general point $p$ of $\ell$ and 
 four lines lying in a general plane $\Pi$ through $\ell$.

\begin{enumerate}

\item[(i)] Suppose that $R$ is a cone with apex $p$ lying on $\ell$, that
    $S_1$ and $S_2$ are inscribed in $R$ and that $\ell\not\in\rho$, so that
    $\ell$ is tangent to neither $S_1$ nor $S_2$.
    (We treat the specialization $\ell\in\rho$ in (iv) below.)
    Then there will be a plane containing $\ell$ and the centres of $S_1$ and
    $S_2$ as well as two distinct planes through $\ell$ that are tangent to
    both $S_1$ and $S_2$.
    Since each line in $\rho$ meets $\ell$ at $p$ and two lie in
    any plane through $\ell$, there will be 4 lines in $\sigma$ meeting a
    general point of $\ell$ and 2 lines in $\sigma$ lying in a general plane
    $\Pi$ through $\ell$, counted with multiplicity.
    These observations remain valid if the configuration becomes degenerate so
    that $S_1$ and $S_2$ are tangent at $p$ and $\ell$ does not lie in their 
    common tangent plane.

    If, for every plane $\Pi$ through $\ell$, the lines $m_\Pi$ and $m'_\Pi$
    of $\sigma$ lying in $\Pi$ meet on $\ell$, then the centres of $S_1$ and
    $S_2$ both lie on $\ell$.
    (To see this, let $\Pi$ be a plane containing $\ell$ and the centres of
    $S_1$ and $S_2$.)
    In this case, $\sigma$ will have three components as detailed in
    Proposition~\ref{T:deg2}: the two complex double
    lines~\eqref{E:special_lines} in $\G$, as well as 
    the family of lines obtained by rotating $m_\Pi$ around $\ell$.
    Any sphere whose centre lies in $\ell$ and is tangent to one
    $m_\Pi$ will be tangent to all the lines in $\sigma$.

    Now suppose that the centres do not lie on $\ell$.
    Then the two tangent planes to $S_1$ and $S_2$ through $\ell$ each contain a
    single line of $\sigma$ with multiplicity 2 (these are also double lines in 
    $\rho$).
    These double tangent lines span a plane containing the centres of
    $S_1$ and $S_2$. 
    There is a unique plane $\Pi$ through $\ell$ such that the two tangent
    lines contained in it intersect on $\ell$---this plane is perpendicular to
    the plane containing $\ell$ and the centres of the two spheres. 
    Thus we have two distinct planes containing the two centres and so their
    intersection gives the line $m$ passing through the two centres. 
    Rotating the known double tangent lines around $m$ we obtain the cone,
    cylinder or double plane parametrized by the omitted degree 2 or degree 1
    component of $\tau$, and thus $\tau$. 
    By Theorem~\ref{T:tau}, this is sufficient to determine $S_1$ and
    $S_2$.\smallskip 

\item[(ii)] Suppose that $\rho$ consists of lines tangent to a smooth conic
    $C=\Pi\cap S_1=\Pi\cap S_2$ in a plane $\Pi$ through $\ell$, and that
    $\ell$ is not tangent to $C$.
    (We treat the specialization of $\ell$ tangent to $C$ in (v) below.)
    Any point $p$ of $\ell$ contains two lines from $\rho$ tangent to the conic,
    and so it contains two lines from $\sigma$.
    We also conclude that any plane through $\ell$ (other than $\Pi$) will
    have four lines from $\sigma$.

    As in the proof of Theorem~\ref{T:tau}, this determines the
    spheres $S_1$ and $S_2$.
    Lastly, these observations remain valid if the positions of $S_1$ and $S_2$
    become degenerate so that they are tangent to  each other at the same point
    of the plane $\Pi$.\smallskip

\item[(iii)] Suppose now that we are in case (iii) of Section~\ref{S:conic}:
    $\ell$ is tangent to the spheres and $\rho$ is a ruling of the hyperboloid
    obtained by revolving $\ell$ about the line joining the centres of the
    spheres. 
    Since $\ell\not\in\rho$, we see that $\ell\in\sigma$.
    Then every point of $\ell$ contains a line in $\rho$ and every
    plane $\Pi$ through $\ell$ also contains a line in $\rho$.
    Thus any point $p$ of $\ell$ lies on three lines in $\sigma$ and any plane
    $\Pi$ through $\ell$ has three lines from $\sigma$.
    For both, the three include $\ell$ itself.

    Suppose that a plane $\Pi$ through $\ell$ is tangent to one of the
    spheres. 
    Then the lines of $\sigma$ in $\Pi$ consist of $\ell$ counted with
    multiplicity 2 and one other tangent line meeting $\ell$ at the point where
    it is tangent to that sphere. 
    In this way we can determine the points at which the spheres are tangent to
    $\ell$. 
    A general plane $\Pi$ through $\ell$ will contain two common tangent
    lines to $S_1$ and $S_2$ apart from $\ell$. 
  {}From these lines and from the known points of the sphere's tangency with
    $\ell$, we can determine the circles $\Pi\cap S_1$, and $\Pi\cap S_2$.
    As before this is sufficient to determine $S_1$ and $S_2$.\smallskip 

\item[(iv)]  The most intricate case is when the envelope $R$ of $\rho$ is a
    cone with apex $p$ tangent to both spheres with $\ell$ a ruling of $R$, so
    that $\ell\in\rho$.
    Thus $\ell$ is tangent to both spheres with common tangent plane $H$. 
    We first observe that $\ell\in\sigma$.  
    To see this, note that for any plane $\Pi$ through $\ell$ (except $H$),
    there will be two lines (besides $\ell$) in $\Pi\cap\Sigma$.
    Thus $H$ contains at least two lines from $\sigma$, counted
    with multiplicity.
    Since $H\cap\Sigma=\ell$ as $H$ is tangent to both
    spheres, we see that $\ell\in\sigma$.
    Thus for $\Pi$ through $\ell$ with $\Pi\neq H$, $\Pi\cap\Sigma$ has two
    lines, and $H$ is the unique plane such that $H\cap\Sigma$ is a
    single line. 
    The plane $M$ through $\ell$ that is perpendicular to the tangent
    plane $H$ contains the centres of the spheres.
    As observed, $M\cap\Sigma$ contains two lines $m$ and $m'$ besides $\ell$.
    The lines $\ell$, $m$ and $m'$ are tangent to the equatorial circles
    $C_1:=M\cap S_1$ and $C_2:=M\cap S_2$ of the spheres.
    Since the fourth line tangent to the circles lies in the cone $R$ with apex $p$,
    the line through their centres (and $p$) meets the lines $m$ and $m'$ at
    their intersection, and in fact bisects {\it one} of the angles there.
    Since there is another pair $C_1'$ and $C_2'$ of circles tangent to the
    three lines with the line connecting their centres bisecting the other
    angle, there are exactly two possibilities for the equatorial circles $C_1$
    and $C_2$, and hence the spheres $S_1$ and $S_2$.
    This geometry in the plane $M$ is indicated in the picture on the left of
    Figure~\ref{Fig:4}. 
\begin{figure}[htb]
\[ 
  \setlength{\unitlength}{1pt}
  \begin{picture}(135,135)(0,-10)
    \put(3.8,0){\epsfxsize=120pt\epsffile{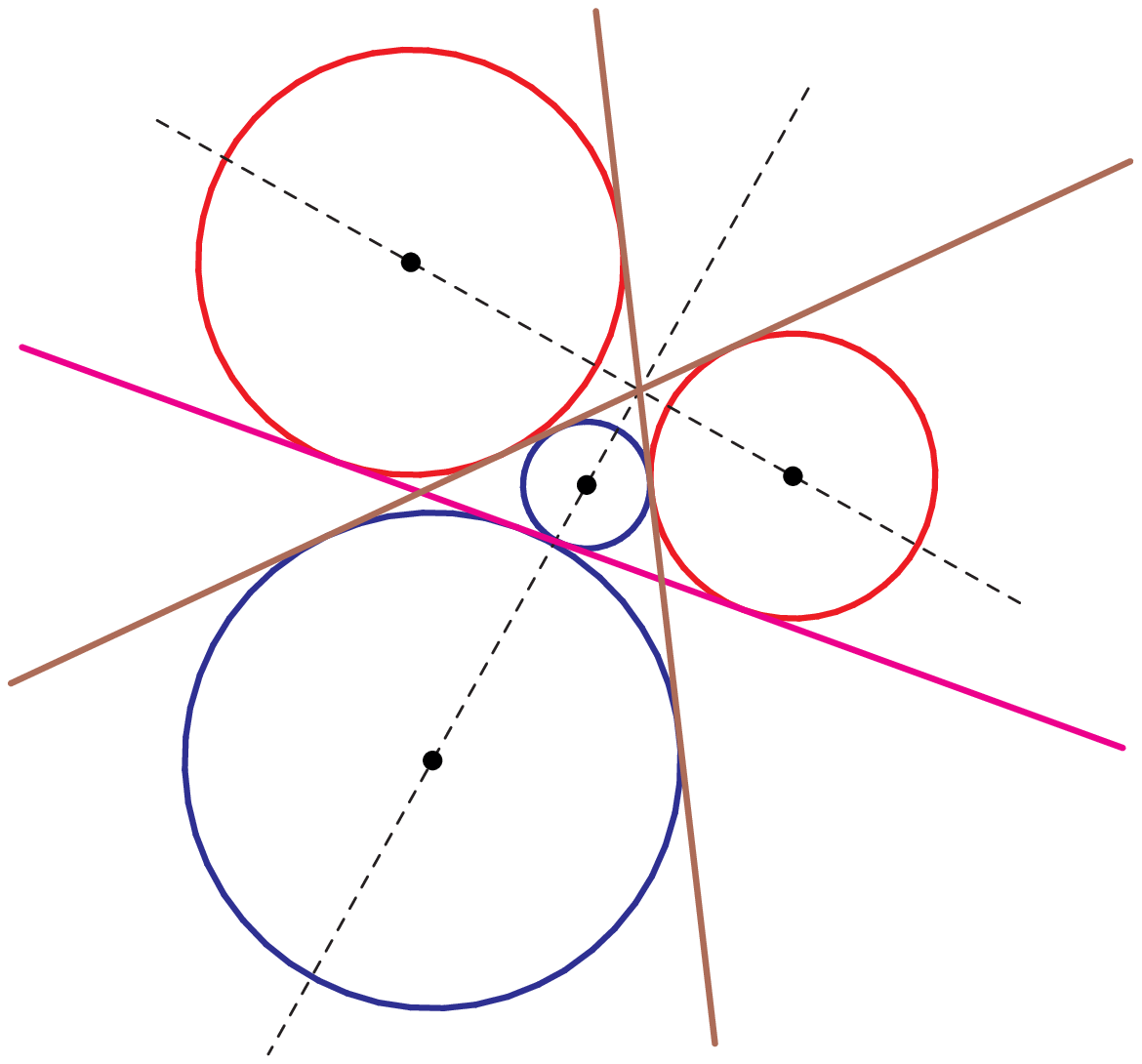}}
    \put( 5,95){$C_1$}  \put(115, 75){$C_2$}
    \put(120,30){$\ell$}
    \put(5,20){$C_1'$}
    \put( 5,64){$C_2'$}\put(20,68){\vector(1,0){42}}
    \put(76,116){$m$}
    \put(115,105){$m'$}
    \put(30,-15){The plane $M$}
  \end{picture}
\qquad\qquad
  \begin{picture}(155,135)(0,-15)
    \put(25,0){\epsfxsize=122.9pt\epsffile{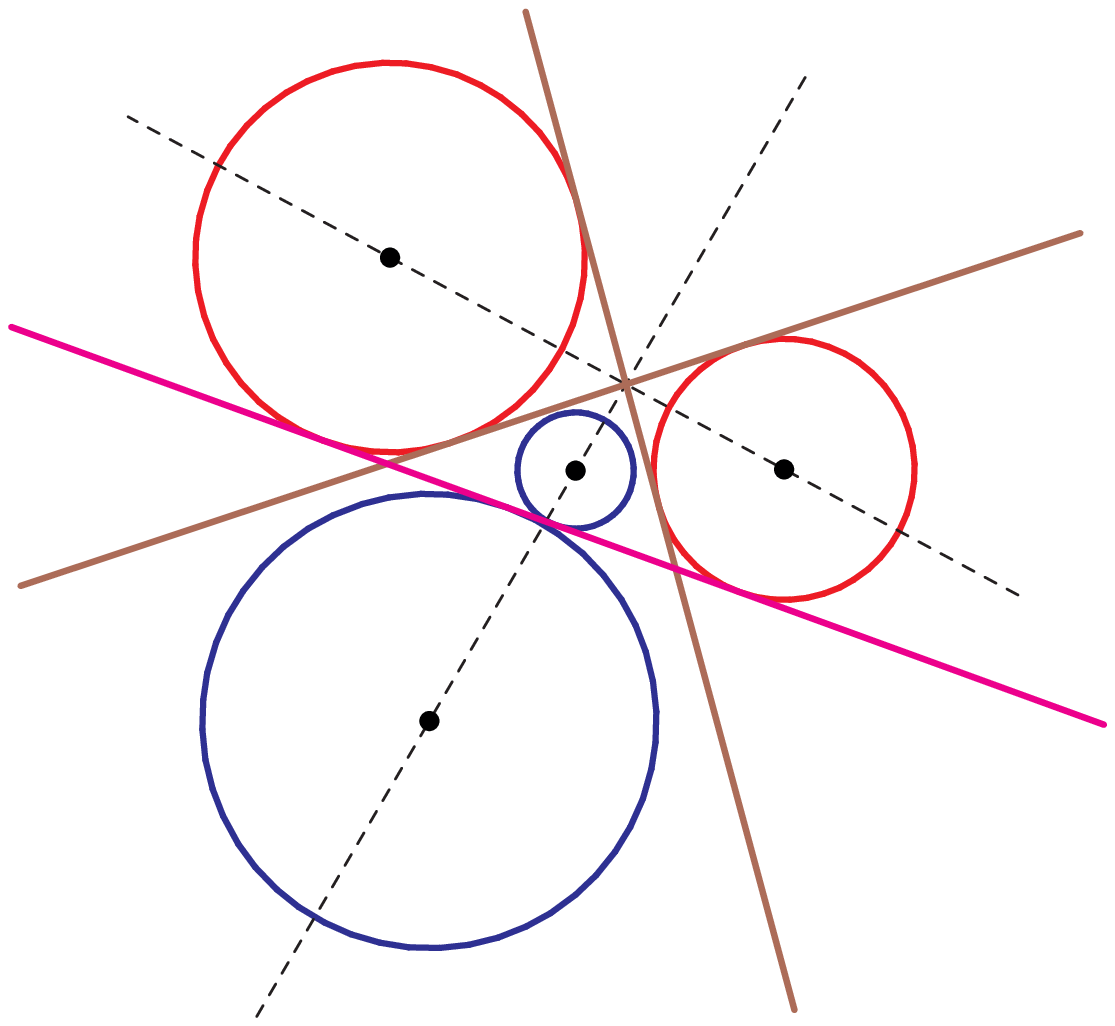}}
    \put( 5,90){$\Pi\cap S_1$}  \put(133, 65){$\Pi\cap S_2$}
    \put(140,20){$\ell$}
    \put( 5,25){$\Pi\cap S_1'$}
    \put( 5,58){$\Pi\cap S_2'$}\put(41,62){\vector(1,0){40.5}}
    \put(90,110){$n$}
    \put(135,95){$n'$}
    \put(40,-20){The plane $\Pi$}
  \end{picture}
\]
\caption{Configurations within planes}
\label{Fig:4}
\end{figure}

The envelope $\sigma$ selects the correct pair.
Indeed, let $S_1,S_2$ and $S_1', S_2'$ be the two possible 
pairs of spheres, and consider another plane $\Pi$ through $\ell$, but
distinct from both $M$ and the tangent plane $H$.
It will intersect these four spheres in four circles, each tangent to $\ell$.
However, only two can be tangent to the other two lines $n,n'$ in $\sigma$ lying
in $\Pi$.
In this way we see that $\sigma$ determines $S_1$ and $S_2$.
The configuration within $\Pi$ is indicated in the picture on the right in
Figure~\ref{Fig:4}, while the configurations of the lines in $M$, the two pairs
of spheres and their intersections with $\Pi$ (drawn in outline) is shown in
Figure~\ref{F:four}.  
Note that the line $n$ is not tangent to the sphere $S'_1$ (the sphere in the
foreground).  
\begin{figure}[htb]
\[
  \begin{picture}(252,140)
    \put(2,0){\epsfxsize=230pt\epsffile{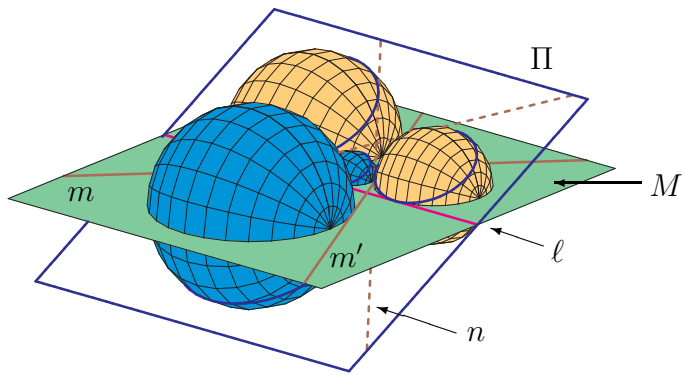}}
    \put(25,65){$m$}   \put(124,40){$m'$}
    \put(205,48){\vector(-3,1){20}}\put(208,42){$\ell$}
    \put(242,72){\vector(-1,0){33}}\put(245,67){$M$}
    \put(200,115){$\Pi$}
    \put(172,15){\vector(-3,1){30}}\put(176,12){$n$}
  \end{picture}
\]
\caption{Configuration of the Spheres}
\label{F:four}
\end{figure}

\item[(v)]  In the last case, $\rho$ consists of lines tangent to a smooth conic
      $C=\Pi\cap S_1=\Pi\cap S_2$ in a plane $\Pi$ through $\ell$, and 
      $\ell$ is tangent to $C$ at a  point $q$.
      We first observe that  $\ell\in\sigma$. 
      To see this, note that for any point $p\neq q$ of $\ell$ there will be two 
      lines in $\sigma$ through $p$ (besides $\ell$).
      Thus $\sigma$ will contain at least two lines that meet $q$, counted with
      multiplicity. 
      Since $\ell$ is the only line tangent to both spheres through the point
      $q$, we conclude that $\ell\in\sigma$.
      Thus $\sigma$ has three lines meeting a point $p\neq q$ of $\ell$
      and $q$ is the unique point of $\ell$ that meets only one line from
      $\sigma$. 

      We complete the proof in this case by noting that every plane $H$
      through $\ell$ (other than $\Pi$) will have 3 lines from $\sigma$.
      This determines the $H\cap S_1$ and $H\cap S_2$ as they are both
      tangent to $\ell$ at the point $q$.
      We conclude once again that $\sigma$ determines the spheres $S_1$ and
      $S_2$, and this completes the proof.\vspace{-16pt}
\end{enumerate}

\end{proof}\medskip

\subsection{The degree of $\sigma$ is 7}

As shown in Section 2.1, any degree 1 component of $\tau$ has multiplicity at
least 2.
Thus there cannot be a component of degree 7.

\section{Quartics}\label{S:Quartics}

We describe completely the configurations of the spheres $S_1$ and $S_2$ in the
case when the curve $\tau$ has an irreducible quartic component $\sigma$.
We also show that there can be no other spheres tangent to all the lines
in $\sigma$.
This uses both algebraic (in complex $\mathbb{P}^3$) and geometric (in
$\mathbb{R}^3$) arguments.

\begin{thm}\label{T:quartics}
 Let $\sigma\subset \G$ be an irreducible real quartic curve
 parametrizing lines meeting a fixed line $\ell$ and
 tangent to a fixed sphere $S_1$. 
 If there is another sphere $S_2$ tangent to all these lines, then either $S_1$
 and $S_2$ are tangent to each other and to $\ell$ at the same point or $S_1$
 and $S_2$ are symmetric with respect to $\ell$ or both.
 In particular, there cannot be a third sphere tangent to all the lines in
 $\sigma$.
\end{thm}

\begin{figure}[htb]
\[
  \epsfxsize=170pt\epsffile{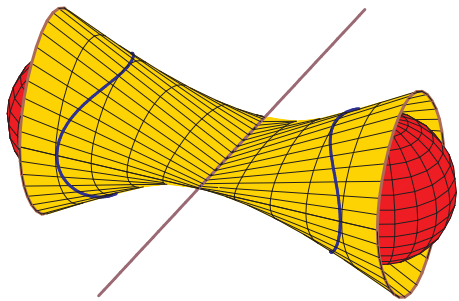}
\]
\caption{Symmetric spheres}
\label{F:symmetric}
\end{figure}

The first case of this theorem is the situation of Example~\ref{E:mult4},
and the second is illustrated in Figure~\ref{F:symmetric}.

For the last statement in Theorem~\ref{T:quartics}, note that if three spheres
are tangent to each other at a point $p$, then one sphere lies inside the other
and the set of common tangents are the lines tangent to the spheres at the point
$p$, which has degree 1 and not 4.
Also we cannot have three distinct spheres that are
pairwise symmetric around $\ell$.

Our proof considers the possibilities for such a degree 4 component of $\tau$,
ruling out most of them.
To begin, there are three possibilities for the dimension of the linear span $H$
of the curve $\sigma$ in Pl\"ucker space, and we treat each separately 
below. 
Note that $H\subset\Lambda_\ell$.

\subsection{The dimension of $H$ is 2}
 In this case, $\sigma$ is a plane quartic.
 Since $\G$ is a quadric, it must contain the plane $H$,
 but the lines lying in such a plane that are tangent to a sphere are
 parametrized by a conic, and not a quartic.
 We conclude that this case is impossible.

\subsection{The dimension of $H$ is 3}\label{S:4.2}
 If $\sigma$ spans a 3-dimensional subspace $H$, then by
 Remark~\ref{R:Three} there are two possibilities: either $H\cap\G$ is
 the set of lines meeting $\ell$ and another (uniquely defined) line
 $\ell'$, or $H$ is tangent to $\G$ at the point $p_\ell$.  In the
 first case, both $\ell$ and $\ell'$ are real and $\sigma$ is a
 quartic curve of lines meeting them that are also tangent to two
 spheres.  This situation was studied in~\cite{2l2s}: by \cite[Theorem
 16]{2l2s} $\ell'$ cannot lie in affine space because those cases do
 not give a quartic curve of common tangents and transversals, and by
 \cite[Theorem 20]{2l2s}, $\ell'$ must be the line at infinity
 perpendicular to $\ell$ and the two spheres must be symmetric with
 respect to $\ell$. This is the second case of
 Theorem~\ref{T:quartics}.  Furthermore, \cite[Theorem 20]{2l2s} also
 implies that there is no other sphere tangent to every line in
 $\sigma$.

 Suppose now that $H$ is tangent is to $\G$ at $p_\ell$.
 Let $\Lambda$ be another hyperplane through $H$ such that
 $H=\Lambda_\ell\cap \Lambda$.
 Then the curve $\sigma$ is the intersection of $\Lambda_\ell$,
 $\Lambda$, $\G$ and either the quadric surface given by $\wedge^2 S_1$ or the
 surface given by $\wedge^2 S_2$. 
 In particular, if $h$ is a linear form defining $H$ in $\Lambda_\ell$, then 
 there is a linear form $k$ such that the quadratic form given by $\wedge^2 S_2$ is a
 linear combination of the quadratic form given by $\wedge^2 S_1$, the Pl\"ucker
 relation~\eqref{E:Pluecker} defining $\G$ and $hk$, modulo $\Lambda_\ell$. 

 Suppose that $\ell$ is the $x$-axis so that $\Lambda_\ell$ is the
 hyperplane $p_{23}=0$. 
 For each $i=1,2$, let $(x_i,y_i,z_i)^T$ be the centre of $S_i$ and let  
 $r_i$ be its radius. 
 Set $q_1$ and  $q_2$ be the quadratic forms given by $\wedge^2 S_1$ and
 $\wedge^2 S_2$ with $p_{23}$ set equal to zero and let $g$ be the Pl\"ucker
 relation with $p_{23}=0$. 
 These define the restrictions of the quadratic forms $\wedge^2 S_1$, 
 $\wedge^2 S_2$ and the Pl\"ucker relation to $\Lambda_\ell$.
 Let $h$ be the equation of $H$ in $\Lambda_\ell$. 
 By the remark at the end of the last paragraph, there exists a real linear form
 $k$ on $\Lambda_\ell$ and real numbers $\lambda, \mu$ and $\nu$
 such that $\lambda q_1+\mu  q_2+\nu g=h k$.
 Furthermore, as either of $q_1$ or $q_2$ may be written in
 terms of the other forms, we have $\lambda,\mu\neq 0$.

 Thus the quadratic form $\lambda q_1+\mu  q_2+\nu g$ has rank 2 and 
 therefore every $3\times 3$ minor of its representation matrix $M$ vanishes. 
 We argue that this is impossible.
 Here is $M$:
\renewcommand{\l}{\lambda}
\newcommand{\m}{\mu}
\newcommand{\n}{\nu}
{\small
 \[ 
  \left(\!\!
   \begin{array}{ccccc}
     \!\begin{array}{r}\l(y_1^2+z_1^2-r_1^2)\\
        +\m(y_2^2+z_2^2-r_2^2)\rule{0pt}{12pt}\end{array} & 
     \!-\l x_1y_1-\m x_2y_2&\!-\l x_1z_1-\m x_2z_2&\l y_1+\m y_2&\l z_1+\m z_2\\
     \rule{0pt}{19pt}
     \!-\l x_1y_1-\m x_2y_2&\!\!\begin{array}{r}\l(x_1^2+z_1^2-r_1^2)
      \\+\m(x_2^2+z_2^2-r_2^2)\rule{0pt}{12pt}\end{array} & 
           \!-\l y_1z_1-\m y_2z_2 &\!-\l x_1-\m x_2 & -\n  \\
     \rule{0pt}{19pt}
     \!-\l x_1z_1-\m x_2z_2& \!-\l y_1z_1-\m y_2z_2 & 
       \!\!\begin{array}{r}\l(x_1^2+y_1^2-r_1^2)\\
       +\m(x_2^2+y_2^2-r_2^2)\rule{0pt}{12pt}\end{array}&\n
      &\!-\l x_1-\m x_2\\
     \rule{0pt}{15pt}
     \l y_1+\m y_2  & -\l x_1-\m x_2 & \n & \l+\m & 0 \\
     \rule{0pt}{15pt}
     \l z_1+\m z_2 &-\n & -\l x_1 -\m x_2 & 0 & \l+\m 
   \end{array}
  \right)
 \]
}
 Since $H$ contains the point
 $p_\ell = (1,0,0,0,0,0)^T$, the linear form $h$ has no $p_{01}$ term, 
 therefore $hk$ has no $p_{01}^2$ term. 
 Thus the upper left entry of $M$ vanishes,
 giving the equation
\begin{equation}\label{E:first}
      \l(y_1^2+z_1^2-r_1^2) +\m(y_2^2+z_2^2-r_2^2)\ =\ 0\,.
\end{equation}
 Let $M_{\text{\it abc},\text{\it def}}$ be the minor formed by the rows 
 {\it abc} and the columns {\it def}. 
 Since $M_{145,145}=(\l+\m)((\l y_1 +\m y_2)^2+(\l z_1 +\m z_2)^2)=0$,
 we have two cases to consider.\medskip

 \noindent{\bf Case 1.} Suppose that $\lambda+\mu=0$.
  Scaling, we may assume that $\lambda=1$ and $\mu=-1$.
  Then $M_{234,345}=\nu(\nu^2+(x_1-x_2)^2)$ and 
  $M_{235,345}=(x_2-x_1)(\nu^2+(x_1-x_2)^2)$.
  Since $x_1$, $x_2$ and $\nu$ are real, this implies that 
  $\nu=0$ and $x_1=x_2$.
  Setting $x:=x_1=x_2$, the matrix $M$ becomes
 \[ 
  \left(\!\!
   \begin{array}{ccccc}
     0 &  - xy_1+ xy_2&- xz_1+ xz_2& y_1- y_2& z_1- z_2 \\
     - xy_1+xy_2&\  z_1^2-z_2^2 + r_2^2-r_1^2 \ & 
             - y_1z_1+y_2z_2 &  0 & 0   \rule{0pt}{14pt} \\
     -xz_1+xz_2& -y_1z_1+ y_2z_2 & 
         y_1^2-y_2^2+ r_2^2-r_1^2&0&0   \rule{0pt}{14pt} \\
      y_1- y_2  & 0 & 0 & 0 & 0    \rule{0pt}{14pt} \\
      z_1- z_2 & 0 & 0 & 0 & 0   \rule{0pt}{14pt} 
   \end{array}
  \right)
 \]
 Since the spheres are not concentric, at least one of $y_1\neq y_2$ or
 $z_1\neq z_2$ holds.
 Suppose that $y_1\neq y_2$.
 Then the equations $0=M_{124,124}=M_{134,134}=M_{124,134}$ imply that the
 entries not in the first row or column vanish, that is
\[
   z_1^2-z_2^2 + r_2^2-r_1^2\ =\ 
     y_1^2-y_2^2+ r_2^2-r_1^2\ =\ y_1z_1-y_2z_2\ =\ 0\,.
\]
 Subtracting the first from~\eqref{E:first} gives
 $0=y_1^2-y_2^2$ and therefore $y_1=-y_2$.
 With the third, we conclude that $z_1=-z_2$ and thus $r_1=r_2$.
 (We reach the same conclusion from $z_1\neq z_2$.)

 Thus the spheres are symmetric with respect to the $x$-axis, 
 the configuration of Figure~\ref{F:symmetric}.
 That envelope corresponds to a curve of degree
 4 whose linear span has the form $\Lambda_\ell\cap\Lambda_{\ell'}$, and so it
 cannot be the curve $\sigma$.
 The other components are two conics as described in the proof of
 Proposition~\ref{T:deg2}: one is the ruling of the
 cone over the spheres with apex the point of symmetry, and the other consists
 of lines in the $xz$-plane tangent to the conic along which the spheres
 meet, and if the spheres are tangent, the components degenerate to a line of
 multiplicity 4. 
 (This is the case $r=1$ in Example~\ref{E:mult4}).
 This gives a contradiction to the existence of $\sigma$ when
 $\lambda+\mu=0$.\medskip 

 \noindent{\bf Case 2.} Suppose that $\lambda+\mu\neq0$, then 
 $\lambda y_1 +\mu y_2=\lambda z_1 +\mu z_2=0$, since these are real
 numbers and the sum of their squares is 0. 
 Since 
 $M_{145,245}=-(\lambda+\mu)^2(\lambda x_1y_1+\mu x_2y_2)$,
 $M_{145,345}=-(\lambda+\mu)^2(\lambda x_1z_1+\mu x_2z_2)$ and
 $M_{245,345}=-(\lambda+\mu)^2(\lambda y_1z_1+\mu y_2z_2)$, we 
 obtain
\[
   \lambda x_1y_1+\mu x_2y_2\ =\ \lambda x_1z_1+\mu x_2z_2\ =\ 
   \lambda y_1z_1+\mu y_2z_2\ =\ 0\,.
\]
Besides scaling $\lambda$ and $\mu$ by a common scalar, the 
only solutions to these five equations are 
\[
  \begin{array}{c}
   x_1=x_2,\  y_1=y_2=0,\  z_1=-\mu,\  z_2=\lambda\makebox[1pt][l]{\,,}\\
   x_1=x_2,\  z_1=z_2=0,\  y_1=-\mu,\  y_2=\lambda\makebox[1pt][l]{\, or}
    \rule{0pt}{13pt}\\
   y_1=y_2=z_1=z_2=0\,.\rule{0pt}{13pt}
 \end{array}
\]
Either of the first two solutions give
$M_{245,245}-M_{345,345}=\pm\lambda\mu(\lambda+\mu)^3$, which leads to
one of the excluded cases $\lambda=0$, $\mu=0$ or $\lambda+\mu=0$. 

The third solution implies that the centres of $S_1$ and $S_2$ lie on
the $x$-axis. 
As in Proposition~\ref{T:deg2} and in Theorem~\ref{T:deg6}, 
$\tau$ will only have components of degrees 1 and 2.
This concludes the proof of the impossibility of a degree 4 component
spanning a 3-plane in Pl\"ucker space that is tangent to the Grassmannian.

\subsection{The dimension of $H$ is 4}

The remaining possibility is that $\sigma$ spans a 4-dimensional subspace.
As $\sigma$ has degree 4, it is necessarily a rational normal quartic
curve.  
There are two possibilities for the geometry of the curve $\gamma\subset
S_1$ along which the lines of $\sigma$ are tangent to $S_1$. 
These lead to restrictions on the possible configurations of 
$S_1$ and $\ell$.
We then consider planes $\Pi$ through $\ell$ 
that are tangent to $S_1$ and contain lines from $\sigma$ (these are determined
from $\sigma$).  
{}From this analysis, we see that the only
possibility for there to be a rational quartic $\sigma$ of lines tangent to two
spheres that also meet a fixed line $\ell$ is when the spheres are tangent to
each other and to $\ell$ at the same point, and thus $\sigma$ is the degree 4
component studied in Example~\ref{E:mult4}.

First, let $\Sigma\subset\mathbb{P}^3$ be the envelope of the rational curve
$\sigma$, a rational ruled surface of degree 4.
The lines in $\sigma$ which are not contained in $S_1$ are tangent to $S_1$ at
well-defined points, and the closure of these points of tangency forms a
curve $\gamma$ lying on $S_1$.
Since $\gamma$ is a real component of multiplicity 2 in the $(4,4)$-curve
$\Sigma\cap S_1$, either $\gamma$ has bidegree
$(1,1)$, so it is a plane conic, or it has bidegree $(2,2)$.
When it has bidegree $(2,2)$, we must have that $\Sigma\cap S_1=\gamma$, in
particular, $\sigma$ contains no lines lying in $S_1$.
When $\gamma$ has bidegree $(1,1)$, $\Sigma\cap S_1$ properly contains $\gamma$,
and the residual curve has bidegree $(2,2)$.
Thus it is the union of lines from the rulings of $S_1$, and so $\sigma$
contains a complex conjugate pair from each ruling of $S_1$. 
The association $\gamma\ni m\cap S_1 \mapsto m\in\sigma$ is the restriction of the
birational map $\Phi\colon S_1\to\G$ studied in Section~\ref{S:birat} to
$\gamma$, so that $\gamma$ is also a rational curve.
Irreducible rational curves of bidegree $(2,2)$ on a quadric 
have a single singularity that is necessarily a simple node or a cusp.
(See the discussion in Section~2 of~\cite{2l2s} and the references therein.)

\begin{lemma}\label{L:config}
 There are two possibilities for the configuration of $S_1$, $\ell$ and
 $\gamma$.
 \begin{enumerate}
  \item[(i)]
            The curve $\gamma$ has bidegree $(1,1)$ and does not contain any
            point where $\ell$ meets $S_1$ or where a plane through $\ell$ is
            tangent to $S_1$.
  \item[(ii)]
            The line $\ell$ is tangent to the sphere $S_1$, 
            the curve $\gamma$ has bidegree $(2,2)$ and is singular at this point of
            tangency.
 \end{enumerate}
\end{lemma}

\begin{proof}
 A point where $\ell$ meets $S_1$ or where a plane through $\ell$ is
 tangent to $S_1$ is a {\it special point}.
 Either $\gamma$ has bidegree $(1,1)$ or
 bidegree $(2,2)$, and either $\ell$ is tangent to $S_1$ or it is not.
 Suppose $\gamma$ has bidegree $(1,1)$ so that it has degree 2.
 By the degree calculation of Lemma~\ref{L:degree},
 $\gamma$ cannot contain any of the special points, proving case (i).

 Suppose that $\gamma$ has bidegree $(2,2)$, so that it has degree 4.
 By Lemma~\ref{L:degree}, $\gamma$ must contain a special point.
 In particular, if $\ell$ is tangent to $S_1$ at the point $p$, then
 $p$ has multiplicity 2 in $\gamma$---that is, $p$ is
 the singular point of $\gamma$.

 We exclude the remaining case of $\gamma$ having bidegree $(2,2)$ and $\ell$
 not tangent to $S_1$.
 By Lemma~\ref{L:degree} (i), either
 $\gamma$ contains all four special points and is nonsingular at each, or
 $\gamma$ contains three of the four special points and is singular at one of
 the three. 
 If $\gamma$ is singular at a special point $q$, then mult$_q\gamma\geq 2$.
 Let $\Pi$ be the tangent plane to $S_1$ at $q$.
 Then the intersection multiplicity of $\Pi$ with $\gamma$ at this point is at
 least 
\[
   2 \cdot {\rm mult}_q \gamma\ \geq \ 4.
\]
 The factor 2 is due to the tangency of $\Pi$ at $\gamma$.
 Since $\gamma$ has degree 4, its total intersection multiplicity 
 with $\Pi$ is 4.
 This implies that $\Pi\cap\gamma =q$.
 Since the tangent plane at any special point contains two others, 
 $\gamma$ cannot be singular at a special point.

 Suppose now that $\gamma$ contains four special points and is thus nonsingular
 at each.
 Then the singular point $q$ of $\gamma$ is not a special point,
 and so the rational map $\Phi$ of~\eqref{E:correspondence} is regular at $q$
 (well-defined on a neighbourhood of $q$).
 It follows that the image of the singular point $q$ of $\gamma$ under $\Phi$ is
 a singular point of the image $\sigma$ of $\gamma$.
 Since any rational normal curve, and in particular $\sigma$ is smooth, this
 contradiction excludes this case. 
\end{proof}

\begin{figure}[htb]  
\[
  \epsfxsize=2.3in\epsffile{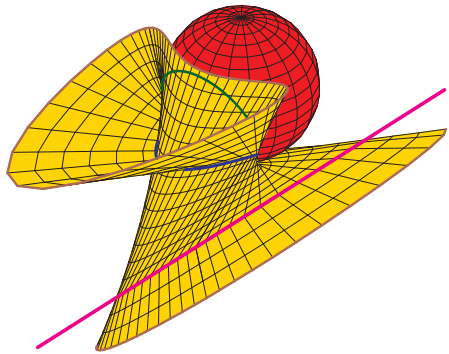}
\]
 \caption{Envelope of a quartic component}\label{F:possibilities}
\end{figure}

Figure~\ref{F:possibilities} displays the envelope of a degree 4 component
whose curve $\gamma$ of tangency satisfies (i) in Lemma~\ref{L:config}.
An envelope corresponding to (ii) is illustrated in
Figure~\ref{F:RatNorm}. 

We complete the proof of Theorem~\ref{T:quartics} by 
considering special lines in $\sigma$ lying in planes $\Pi$ through $\ell$.
This allows us to distinguish the two cases for the curve $\gamma$ given in 
Lemma~\ref{L:config}, and then to show that the only possibility for two spheres 
to be tangent to the lines of $\sigma$ is for the spheres to be tangent to each
other and to $\ell$ at the same point, the configuration of
Example~\ref{E:mult4}.  

We first review some complex Euclidean analytic geometry.
Suppose that $\ell$ is the $x$-axis.
For $\lambda,\mu\in\mathbb{C}$ satisfying $\lambda^2+\mu^2=1$, consider the 
plane through $\ell$ given by the equation $\lambda y=\mu z$.
All planes through $\ell$ have this form, except for the two planes 
$y=\pm iz$, and these two are independent of coordinate choice.
Euclidean coordinates for the plane $\Pi\colon\lambda y=\mu z$ are given by $x$
and $v:=\mu y + \lambda z$.
Suppose that $S_1$ is a sphere with centre $(x_0,y_0,z_0)$ and radius $r$.
Then $\Pi\cap S_1$ is the circle
 \begin{equation}\label{E:circle}
   (x-x_0)^2 + (v-v_0)^2\ =\ r^2-y_0^2-z_0^2+v_0^2\,,
 \end{equation}
where $v_0=\mu y_0+\lambda z_0$.
The plane $\Pi$ is tangent to $S_1$ when $v_0^2=y_0^2+z_0^2-r^2$, with the
circle degenerating to the two lines $(x-x_0)=\pm i(v-v_0)$ with slopes 
$\pm i$.
We remark that a line in such a plane with slope $\pm i$ is tangent to a circle
only when the circle is degenerate.

Consider the planes $\Pi$ through $\ell$ having an equation
of the form $\lambda y=\mu z$ with $\lambda^2+\mu^2=1$ that contain a line from 
$\sigma$ with slope $i$ or $-i$.
Call such a plane {\it distinguished}.
Note that the distinguished planes are determined by $\sigma$.
Since distinguished planes are tangent to the sphere $S_1$, there are at
most two. 
We  consider the three cases of two, one or zero distinguished planes
through $\ell$.
Recall that if $\gamma$ has bidegree $(2,2)$, then $\Sigma\cap S_1$ necessarily
equals $\gamma$ and there are no lines from $\sigma$ lying in $S_1$.
Thus there can be no distinguished planes if $\gamma$ has bidegree $(2,2)$.\smallskip

\noindent{\bf Two distinguished planes. }
We must be in case (i) of Lemma~\ref{L:config}.
Furthermore, the line $\ell$ does not contain the centre of $S_1$, for then no
plane through $\ell$ of the form $\lambda y=\mu z$ with $\lambda^2+\mu^2=1$ is
tangent to $S_1$ (tangent planes to $S_1$ through a line containing its centre
have the form $y=\pm iz$). 
Conversely, if we are in case (i) of Lemma~\ref{L:config} and $\ell$ is neither
tangent to $S_1$ nor contains the centre of $S_1$, then any line in $\sigma$ lying
in a tangent plane $\Pi$ to $S_1$ through $\ell$ must have slope $\pm i$:
otherwise such a line in $\sigma$ meets $S_1$ only at the point $p$ of tangency
of $\Pi$ to $S_1$, and thus $p\in\gamma$, which contradicts (i) of
Lemma~\ref{L:config}. 

Consider lines in $\sigma$ lying in a distinguished plane $\Pi$.
If there is only a single such line $m$, then its complex conjugate
line $\overline{m}$ also lies in $S_1$ and is also in $\sigma$.
Necessarily $\overline{m}\subset\overline{\Pi}$, the plane complex conjugate to
$\Pi$.
Since $m$ is the only such line in $\Pi$, we have that $\Pi\neq\overline{\Pi}$.
As $m$ and $\overline{m}$ are complex conjugate, they lie in different
rulings of $S_1$, and therefore meet in a point $p$ (necessarily real) of $S_1$.
Note that $p$ lies in the intersection of $\Pi$ with $\overline{\Pi}$, which is
$\ell$. 
Furthermore, the plane $H$ spanned by $m$ and $\overline{m}$ is tangent to
$S_1$ at $p$.

If we have a second sphere tangent to all the lines of $\sigma$, then
it is also tangent to $H$ at $p$.
However, this configuration the spheres and line $\ell$ is exactly that of
Theorem~\ref{T:deg1} (i) (see also (ia) in Figure~\ref{F:Thm1}). 
In this case, the lines in $H$ through $p$ form a component of degree 1 and
multiplicity 2 in the curve $\tau$ of tangents to the two spheres that meet
$\ell$. 
By Proposition~\ref{T:deg2}, the residual sextic is irreducible as
$\ell$ neither contains the centres of the spheres nor are they symmetric about
$\ell$. 
However, this contradicts the existence of the quartic component $\sigma$.
Observe that this argument only depends upon there being exactly one line from 
$\sigma$ in one of two tangent planes to $S_1$ through $\ell$.

Now suppose that $\sigma$ contains two lines in each distinguished plane $\Pi$.
Since the two lines from $\sigma$ lying in a distinguished plane meet at the point of
tangency, $\sigma$ determines two planes through $\ell$ tangent to $S_1$ and the
points of tangency. 
These data determine the sphere $S_1$ and there is no possibility for a second
sphere.\smallskip 

\noindent{\bf One distinguished plane. }
Now suppose that there is a unique distinguished plane $\Pi$.
Then $\gamma$ has bidegree $(1,1)$ and $\ell$ does not contain the centre of
$S_1$. 
Since $\sigma$ must contain lines from each ruling of $S_1$ and these must lie
in a plane tangent to $S_1$ through $\ell$, we conclude that $\sigma$ contains
both lines in $\Pi$ lying in $S_1$.
Furthermore, $\Pi$ is real as otherwise $\overline{\Pi}$ is a different
distinguished plane.
Since the two lines from $\sigma$ meet at the point $p$ of tangency of $S_1$ to
$\Pi$, this point and $\Pi$ are determined by $\sigma$.

If there is a second sphere tangent to all the lines from $\sigma$,
it is necessarily tangent to $\Pi$ and hence also to $S_1$ at the point $p$.
If $\ell$ is not tangent to $S_1$, this is the configuration of
Theorem~\ref{T:deg1} (ii) (see also (ib) in Figure~\ref{F:Thm1}),
and, as the spheres are not symmetric about $\ell$, Proposition~\ref{T:deg2}
leads to a contradiction as before.
On the other hand $\ell$ cannot be tangent to $S_1$, for then this is the
configuration of Example~\ref{E:mult4}, and in that case lines lying in the
common tangent plane had real slopes whereas the lines here have imaginary
slopes.
\smallskip 

\noindent{\bf No distinguished planes. }
Suppose that no plane through $\ell$ of the form $\lambda y=\mu z$ with
$\lambda^2+\mu^2=1$ contains a line from $\sigma$ with slope $i$ or
$-i$. 

If the curve $\gamma$ has bidegree $(1,1)$, then $\ell$ necessarily
contains the centre of $S_1$, as $\sigma$ contains lines lying in $S_1$ and
these must lie in a plane through $\ell$ that is tangent to $S_1$.
Otherwise $\ell$ is tangent to $S_1$ at the singular point of the $(2,2)$-curve  
$\gamma$.  
We distinguish these two cases by considering a plane $\Pi\colon y=\pm iz$
through $\ell$,
which is independent of the choice of coordinates (up to complex conjugation).
If the centre of $S_1$ lies on $\ell$, then $\Pi\cap S_1$ will be two lines
meeting at the point $(0,0,1,\pm i)^T$ at infinity.
Since $\gamma$ does not contain the point of tangency, at least one of these
lines lies in $\sigma$.
As argued in the case of two distinguished planes, we cannot have $\sigma$
containing only one of these lines, so $\gamma$ must contain both lines,
and they meet at this point.
On the other hand, if $\ell$ is tangent to $S_1$, then $\Pi\cap S_1$ will be a
smooth conic containing the point $(0,0,1,\pm i)^T$ at infinity, and so the
tangents in $\sigma$ lying in $\Pi$ cannot meet at this point.
Thus we may distinguish the two cases when there are no distinguished planes.

Suppose that $\gamma$ has bidegree $(1,1)$ so that the lines from $\sigma$ in
the plane $y=\pm iz$ contain the point $(0,0,1,\pm i)^T$.
Then the line $\ell$ necessarily contains the centre of $S_1$. 
Thus if there are two spheres tangent to all the lines from $\sigma$, then their
centres must lie on $\ell$.  
By Proposition~\ref{T:deg2}, the curve $\tau$ of tangents to the two spheres 
that also meet $\ell$ will then have four components, none of which has degree
more than 2, so this case does not occur.

We finally conclude that $\ell$ is tangent to the sphere $S_1$ and $\gamma$ has
bidegree $(2,2)$. 
Since this determination depends only upon $\sigma$, we conclude that the other
sphere $S_2$ is also tangent to $\ell$.
The points of tangency can be different or coincide, and the planes of tangency
can be different or coincide.
If the spheres are tangent to $\ell$ at the same point with the same tangent
plane, then we have the configuration of Example~\ref{E:mult4}.
Thus the curve $\tau$ of common tangents to the spheres that meet $\ell$ consists
of $\sigma$ together with the lines in the common tangent plane through the
point of tangency, a component of multiplicity 4.
By Theorem~\ref{T:tau}, this determines $S_1$ and $S_2$; in particular, there
cannot be a third sphere.

We complete the proof of Theorem~\ref{T:quartics} and thus of
Theorem~\ref{T:ls2}  by disposing of the three remaining possibilities.
If the spheres are tangent to $\ell$ at different points, but with the same
tangent plane, then we have the configuration of (iv) in the proof of
Theorem~\ref{T:deg6}. 
If the points of tangency coincide, but the tangent planes are distinct, this is
the configuration of (v) in the proof of Theorem~\ref{T:deg6}. 
If the points of tangency and the tangent planes are all distinct, then this is
the configuration of Section~\ref{S:conic} (iii).
In each of these cases the curve $\tau$ of common tangents to the spheres
that meet $\ell$ has a component of degree 2, and, by Proposition~\ref{T:deg2},
the rest of $\tau$ consists of an irreducible sextic, contradicting our
assumption of a quartic component $\sigma$.

\def\cprime{$'$}
\providecommand{\bysame}{\leavevmode\hbox to3em{\hrulefill}\thinspace}
\providecommand{\MR}{\relax\ifhmode\unskip\space\fi MR }
\providecommand{\MRhref}[2]{%
  \href{http://www.ams.org/mathscinet-getitem?mr=#1}{#2}
}
\providecommand{\href}[2]{#2}

\end{document}